\documentclass[11pt,dvipdfmx]{amsart}
\usepackage{amscd}
\usepackage{amsfonts}
\usepackage{amssymb}
\usepackage{amsmath}
\usepackage{color}
\usepackage{latexsym}

\numberwithin{equation}{section}

\newtheorem{pro}{Proposition}[section]
\newtheorem{lemma}[pro]{Lemma}
\newtheorem{theorem}[pro]{Theorem}
\newtheorem{definition}[pro]{Definition}
\newtheorem{remark}[pro]{Remark}
\newtheorem{note}[pro]{Note}
\newtheorem{eg}[pro]{Example}

\newtheorem{corollary}[pro]{Corollary}

\def\mbi#1{\boldsymbol{#1}} 

\newcommand{\sfD}{{\mathsf{D}}}
\newcommand{\sfE}{{\mathsf{E}}}

\newcommand{\sfR}{{\mathsf{R}}}

\newcommand{\sfJ}{{\mathsf{J}}}

\newcommand{\Ll}{\mathcal{L}}
\newcommand{\al}{\alpha}
\newcommand{\be}{\beta}
\newcommand{\ga}{\gamma}

\newcommand{\e}{\varepsilon}
\newcommand{\se}{\theta}
\newcommand{\om}{\omega}
\newcommand{\Om}{\Omega}
\newcommand{\we}{\wedge}
\newcommand{\lam}{\lambda}
\newcommand{\sig}{\sigma}
\newcommand{\ra}{{\rightarrow}}
\newcommand{\lra}{{\longrightarrow}}

\newcommand{\ie}{{i.e$.$\,}}
\newcommand{\cf}{{cf$.$\,}}

\newcommand{\CC}{{\mathbb C}}

\newcommand{\HH}{{\mathbb H}}

\newcommand{\NN}{{\mathbb N}}
\newcommand{\PP}{{\mathbb P}}
\newcommand{\RR}{{\mathbb R}}

\newcommand{\ZZ}{{\mathbb Z}}

\newcommand{\cK}{{\mathcal{K}}}
\newcommand{\cR}{{\mathcal{R}}}
\newcommand{\cN}{{\mathcal{N}}}
\newcommand{\cM}{{\mathcal{M}}}
\newcommand{\cA}{{\mathcal{A}}}
\newcommand{\cB}{{\mathcal{B}}}
\newcommand{\cG}{{\mathcal{G}}}
\newcommand{\cV}{{\mathcal{V}}}

\newcommand{\cL}{{\mathcal{L}}}
\newcommand{\cX}{{\mathcal{X}}}

\newcommand{\GL}{\mathop{\rm GL}\nolimits}
\newcommand{\SL}{\mathop{\rm SL}\nolimits}

\newcommand{\PSL}{\mathop{\rm PSL}\nolimits}
\newcommand{\SO}{\mathop{\rm SO}\nolimits}
\newcommand{\Sp}{\mathop{\rm Sp}\nolimits}
\newcommand{\PSp}{\mathop{\rm PSp}\nolimits}

\newcommand{\PO}{\mathop{\rm PO}\nolimits}
\newcommand{\PU}{\mathop{\rm PU}\nolimits}
\newcommand{\U}{\mathop{\rm U}\nolimits}

\newcommand{\Sim}{\mathop{\rm Sim}\nolimits}
\newcommand{\Aff}{\mathop{\rm Aff}\nolimits}
\newcommand{\Aut}{\mathop{\rm Aut}\nolimits}
\newcommand{\E}{\mathop{\rm E}\nolimits}

\newcommand{\dev}{\mathop{\rm dev}\nolimits}
\newcommand{\Dev}{\mathop{\rm Dev}\nolimits}

\newcommand{\Psh}{\mathop{\rm Psh}\nolimits}
\newcommand{\PSO}{\mathop{\rm PSO}\nolimits}
\newcommand{\Diff}{\mathop{\rm Diff}\nolimits}
\newcommand{\Isom}{\mathop{\rm Isom}\nolimits}

\newtheorem*{thm*}{Theoem}
\newtheorem*{thmA}{Theorem A}
\newtheorem*{thmB}{Theorem B}
\newtheorem*{thmC}{Theorem C}
\newtheorem*{thmD}{Theorem D}
\newtheorem*{thmE}{Theorem E}
\newtheorem*{thmF}{Theorem F}

\newtheorem*{thmI}{Theorem I}

\newtheorem*{propA}{Proposition A}

\newtheorem*{propC}{Proposition C}

\newtheorem*{coroG}{Corollary 7}

\begin{document}

\setcounter{tocdepth}{2}
\setcounter{page}{1}

\baselineskip 13pt \pagestyle{myheadings} 

\title[\ ]
{Quaternionic contact structure
with integrable complementary distribution }
\author[]{Yoshinobu KAMISHIMA}
\address{Department of Mathematics, Josai University\\
Keyaki-dai 1-1, Sakado, Saitama 350-0295, Japan}
\email{kami@tmu.ac.jp}



\begin{abstract}
We study positive definite quaternionic contact $(4n+3)$-manifolds 
($qc$-manifold for short).  
Just like the $CR$-structure contains the class of Sasaki manifolds,
the $qc$-structure admits a class of
$3$-Sasaki manifolds with integrable distribution isomorphic to $\mathfrak{su}(2)$.  
A big difference concerning the integrable complementary $qc$-distribution 
$V$ of the $qc$-structure from $3$-Sasaki structure 
is the existence of Lie algebra not isomorphic to $\mathfrak{su}(2)$. 
We take up \emph{non-compact} $qc$-manifolds to find out a
salient feature of topology and geometry in case
 $V$ generates the $qc$-transformations $\RR^3$.
\end{abstract}

\maketitle

\tableofcontents

\thispagestyle{empty}

\section{Introduction}\label{sec:in}
It is known that the quaternionic contact structure ($qc$-structure for short) 
on a $4n+3$-dimensional manifold $(M,\sfD)$
contains a class of $3$-Sasaki manifolds
as a quaternionic $CR$-structure.
In this case the complementary integrable distribution $V$ is locally 
isomorphic to $\mathfrak {su}(2)$ ($\cong$ $\mathfrak{sp}(1)$).  
In order to obtain the Reeb field as in $CR$-structure,
the $qc$-automorphism group $\Aut_{qc}(M)$
is relatively too big to act properly on $M$. Like as the 
pseudo-Hermitian transformation subgroup extracted from 
the $CR$-transformation group of a $CR$-manifold,
we may find a subgroup $\Psh_{qc}(M)$ from $\Aut_{qc}(M)$.
$V$ is said to the \emph{$qck$-distribution}
if it generates a $3$-dimensional subgroup of $\Psh_{qc}(M)$.
The idea of reduction to this group was performed in our paper
\cite{OK1} by using the vanishing of equivariant smooth
cohomology groups. Geometrically this means,
given an ${\rm Im}\, \HH$-valued
$1$-form $\om$ representing the $qc$-structure on $M$,
there is a smooth function $v\in C^\infty(M,\RR^+)$ such that
the one-form $\eta=v\cdot \om$ gives rise to 
the subgroup $\Psh_{qc}(M,\eta)\leq \Aut_{qc}(M)$.
With the aid of the work \cite{RS}, the following theorem is obtained 
along the same method of \cite{OK1} which
clarifies how 
$\Psh_{qc}(M)$ interacts the $qc$-structure.
(See Theorem \ref{qvanish}.)  

\begin{thmA}
Let $(M,(\sfD,\{J_\al\}_{\al=1}^3))$ be a 
positive definite $qc$-manifold.
Then either one of the following holds: 
\begin{enumerate}
\item[(i)] There exists a positive definite
$qc$-structure $(\eta,\{J_\al\}_{\al=1}^3)$\, with $D= \ker\, \eta $\, 
such that
\[\Aut_{qc}(M,\sfD,\,\{J_\al\}_{\al=1}^3\,)=\Psh_{qc}(M,\eta,\,\{J_\al\}_{\al=1}^3\,) \;.\]
\item[(ii)]
 $M$ has a spherical $qc$-structure isomorphic to
either the standard sphere $S^{4n+3}$ or the quaternionic Heisenberg Lie group $\cM$.
\end{enumerate} In case {\rm (ii)}, 
$\displaystyle (\Psh_{qc}(M),\,\Aut_{qc}(M),\, M)$ is exactly the following: 
\begin{equation*}\begin{split}
\begin{cases}
\left(\Sp(n+1)\cdot \Sp(1),\ {\rm PSp}(n+1,1),\ S^{4n+3}\right)\,,\, \\
\left(\cM\rtimes {\rm Sp}(n)\cdot \Sp(1),\
 \cM\rtimes ({\rm Sp}(n)\cdot \Sp(1)\times \RR^+),\ \cM\,\right).\,\\
\end{cases}
\end{split}\end{equation*}
\end{thmA}
A $qc$-manifold $(X,\sfD)$ with $\ker\,\om=\sfD$
is \emph{positive definite} if the Levi form 
$d\om_\al\circ J_\al: \sfD\times \sfD\lra\, \RR$ defined by
$d\om_\al\circ J_\al(\mbi{x},\mbi{y})=d\om_\al(J_\al\mbi{x},\mbi{y})$ is a 
positive definite symmetric bilinear form $(\al=1,2,3)$.
A $qc$-manifold may be assumed to be positive definite throughout this paper. 

A $qc$-manifold $X$ with the $qck$-distribution $V$
generalizes the notion of $3$-Sasaki (quaternionic $CR$-) manifold
(with the Killing Reeb fields), that is $V$ generates $\cR\leq \Psh_{qc}(X)$,
called the \emph{$qck$-group}.
Then  
$$\Psh_{qc}(X)=N_{\Psh_{qc}(X)}(\cR)$$
where $\displaystyle N_{\Psh_{qc}(X)}(\cR)$
is the normalizer of $\cR$ in $\Psh_{qc}(X)$. 
(See Corollary \ref{normalizer}, also \cite{OK}.)
In case $M$ is a compact $qc$-manifold with a $qck$-group
$T\leq\Psh_{qc}(M)$, $T$ is locally isomorphic to $\Sp(1)$ or
$T$ is isomorphic to the toral group $T^3$ (\cf Proposition \ref{class3}).
For $qCR$-manifolds ($3$-Sasaki manifolds),
see \cite{BG1}, \cite{IK}. 
In general see \cite{IMV} for compact $qc$-manifolds and 
\cite{Biq}, \cite{IMV1}, \cite{IV} for the work on the $qc$-structure
and the references therein. 

In this paper we study 
non-compact positive definite $qc$-manifolds $X$ with a non-compact $qck$-group $\cR$,
mainly $\cR=\RR^3$. (Compare Proposition \ref{soR3}.)
We discuss Riemannian submersions obtained from $X$.
Compare Theorem \ref{hcr} (also Proposition \ref{ccstruc}).

\begin{thmB}\label{hcr}
Let $X$ be a simply connected 
non-compact $qc$-manifold with an ${\rm Im}\, \HH$-valued $1$-form $\eta$ and
the $qck$-distribution $V$
generating $\RR^3\leq \Psh_{qc}(X)$.
Then there is a principal Riemannian 
submersion $\displaystyle
\RR^3\ra\ X\stackrel{p}\lra\, Y$ where $(Y,\Om)$ is a simply connected hyperK\"ahler 
manifold such that $p^*\Om=d\eta$.
For each $\RR\leq \RR^3$, the quotient
$X/\RR$ is a \emph{complex contact} manifold which
admits a holomorphic principal bundle over 
a hyperK\"ahler manifold $Y=X/\RR^3$:
\begin{equation}\label{1holbund}
\begin{CD} \CC@>>> X/\RR@>>> Y.
\end{CD}\end{equation}  
For each $\RR^2\leq \RR^3$,
the quotient $X/\RR^2$ is
a strictly pseudoconvex $CR$-manifold admitting
a pseudo-Hermitian $($Sasaki$)$ bundle: 
\begin{equation}\label{PHSA}
\begin{CD}
\RR@>>> X/\RR^2@>>> Y.
\end{CD}
\end{equation}
\vskip0.1cm
Similarly suppose $M$ is a compact $qc$-manifold $X/\Gamma$ 
whose $qck$-group $T^3$ lifts to an $\RR^3$-action to $X$. Then\\
{\rm (i)}\, For each $S^1\leq T^3$,
the quotient oribifold $M/S^1$
supports a complex contact structure such that 
$M/S^1$ is the holomorphic orbibundle
over the hyperK\"ahler orbifold $Z=M/T^3$: 
\begin{equation}\label{holbund1}
\begin{CD}
 T^1_\CC@>>> M/S^1@>>> Z.
\end{CD}\end{equation}
{\rm (ii)}\, For a torus $T^2\leq T^3$,
the quotient orbifold $M/T^2$
admits a strictly pseudoconvex $CR$-structure. Furthermore
this gives the pseudo-Hermitian $($Sasaki\,$)$ orbibundle: 
\begin{equation}\label{holbund2}\begin{CD}
S^1@>>> M/T^2@>>> Z.
\end{CD}\end{equation}
{\rm (iii)}\, The fundamental group $\Gamma=\pi_1(M)$
to \eqref{holbund1} 
induces a nontrivial group extension:
$\displaystyle 1\ra\, \ZZ^2\ra\, \Gamma/\ZZ \lra\, Q\ra 1$ where $Y/Q=Z$.
$\Gamma$ also assigns to \eqref{holbund2} 
a \emph{nontrivial} group extension:
$\displaystyle 1\ra\, \ZZ\ra\, \Gamma/\ZZ^2 \lra\, Q\ra 1$.
Then \eqref{holbund1}
is a nontrivial $T^1_\CC$-orbibundle.
\end{thmB}

\noindent 
In case $X$ is a simply connected $qc$-manifold
with the $qck$-group $\RR^3$, 
the $qc$-Hermitian group (respectively
hyperK\"ahler group) can be described exactly  
as $\displaystyle \Psh_{qc}(X, \om,\{J_\al\}_{\al=1}^3)=\bigl\{
f\in\Diff(X)\, \big|\  f^*\om
=a \cdot\om\cdot \bar a\,, \ a\in\Sp(1)\bigl\}$,
$\displaystyle \Isom_{hK}(Y, \Om,\{{\sf J}_\al\}_{\al=1}^3)=\bigl\{
h\in\Diff(Y)\, \big|\  h^*\Om =b \cdot\Om\cdot \bar b\,,\ b\in \Sp(1)\big\}$.
In general note that $a,\,b$ are smooth functions on $X$, $Y$ respectively.
(Compare Corollary \ref{abeliabcase} for the precise description.)
The following concerns
the structure of $\displaystyle \Psh_{qc}(X)$ (\cf Corollary \ref{exact3}).

\begin{propC}\label{exact3} Let $X$ be a simply connected $qc$-manifold
with $qck$-group $\RR^3$.
There is a natural exact sequence:
$$\begin{CD} 1@>>>\RR^3@>>>\Psh_{qc}(X)@>{\phi}>> \Isom_{hK}(Y)@>>> 1. \end{CD}$$
\end{propC}
\noindent For example if $X=\cM$ is a quaternionic Heisenberg 
Lie group as a $qc$-homogeneous manifold,
then it follows 
$\Psh_{qc}(X)=\cM\rtimes \Sp(n)\cdot \Sp(1)$ where $Y$ is
the quaternionic space $\HH^n$
such that $\Isom_{hK}(Y)=\HH^n\rtimes \Sp(n)\cdot \Sp(1)$.

Given a $qCR$-manifold $X$, the product
$\RR^+\times X$ with the cone metric
is known to admit an (imcomplete) hyperK\"ahler structure (\cf \cite{AK1} for instance).
This construction cannot be applied to $qc$-manifolds $X$ with 
the $qck$-group $\RR^3$ (see Note \ref{actionSpM}). 
However we shall construct a hyperK\"ahler metric on $\RR\times X$
suitable for any $qc$-manifold $X$ with $qck$-group $\RR^3$.
(Compare Proposition \ref{hyperKX}.) Applying this construction to
the quaternionic Heisenberg Lie group $\cM$, we obtain 

\begin{thmD}\label{hyperM}
There is a complete hyperK\"ahler metric $g_0$ on the quaternionic space
$\HH^{n+1}$ $(n\geq 1)$ 
such that the quaternionic isometry group is
$\Isom_{hK}(\HH^{n+1},g_0)=\HH\rtimes (\Sp(n)\times \Sp(1))$. 
In particular, $g_0$ is not equivalent to the standard quaternionic metric
of $\HH^{n+1}$ up to a quaternionic isometry. 
\end{thmD}
Our next purpose is to classify \emph{spherical homogeneous $qc$}-manifolds. 

\begin{definition}\label{def:spuni}
A $4n+3$-dimensional positive definite $qc$-manifold $M$ is spherical
$($or uniformizable$)$
if it is locally modeled on $S^{4n+3}$
with coordinate changes lying in $\PSp(n+1,1)$.
Equivalently there exists a $qc$-developing map of the universal covering $\tilde M$
 to $S^{4n+3}$.
\end{definition}\noindent Compare \cite{IV}.
The classification of homogeneous $qc$-manifolds is
a difficult subject in its own right.
The spherical homogeneous $qc$-manifolds 
 are comparatively nontrivial examples appropriate to Proposition C.    
(Refer to \cite{BS} for the spherical $CR$ case.) See Theorem \ref{isoqchomo}
(\cf Propositions $\ref{tran2}$, $\ref{homGk}$).

\begin{thmE}
Any simply connected spherical homogeneous $qc$-manifold $M$ is $qc$-isomorphic to
$S^{4n+3}$, $S^{4n+3}-S^{4m-1}$\  $(1\leq m\leq n)$ 
or $\cM$, $X(k)$, $X(0)$\  $(1\leq k\leq n)$.   
In particular each one of $\cM$, $X(k)$ or $X(0)$ is a principal $qc$-bundle
 over a domain of $\HH^n$ with $qck$-fiber $\RR^3$.
\end{thmE}
It is important to verify the geometric properties of $G$-structure on $qc$-manifolds
(\cf \cite{Biq}). 
It is worthwhile pointing out 
the vanishing theorem of Biquard's quaternionic conformal geometry
obtained by S. Ivanov and D. Vassilev \cite{IV}.
They proved the vanishing of a $qc$-conformal (curvature) tensor implies 
a $qc$-manifold is spherical, which eventually arrives at Theorem A 
(followed by the works of \cite{RS}, \cite{CF}).
 Associated with the canonical Riemannian metric $g$ of \eqref{Riemann},
we simply calculate the curvatures of 
positive definite $qc$-manifolds with $qck$-group $\RR^3$ 
by applying O'Neill's formula to the fiberings of Theorem B. 
Then $X$ is no longer an Einstein manifold unlike a $3$-Sasaki manifold.
Compare Propositions \ref{Drist}, \ref{notEinstein}.  
Finally we prove the geometric uniqueness of $qc$-manifolds $X$
with $qck$-group $\RR^3$ concerning the model quaternionic 
Heisenberg Lie group $\cM$ under several conditions.
 This assertion is brought together  
 by Propositions \ref{Dristen}, \ref{conthomoqc}, 
\ref{homoqc1}, \ref{compactuni}.
(See also \cite[Theorem 4.4]{KA2}.)

\vskip0.1cm
\emph{The paper is organized as follows}.
In Section \ref{QContact group}, we prepare several basic facts of $qc$-structures,
especially the equivalence relation of $qc$-conformal change.
Section \ref{Confinv} is a review of the $qc$-conformal invariant
(\cf \cite{OK}) whose vanishing on a $qc$-manifold $M$ reduces the $qc$-automorphism group 
$\Aut_{qc}(M,\om)$ to 
the $qc$-Hermitian group $\Psh_{qc}(M,\eta)$. 
In Section \ref{complemintg} we introduce a $qck$-distribution on $qc$-manifolds.
In general there is a complementary distribution $V$ to
the $qc$-subbundle $\sfD$ of a $qc$-manifold $M$. We study
a $qc$-manifold whose distribution generates 
a $qck$-group $\cR$ in $\Psh_{qc}(M)$.
In Section \ref{sec:fib}, we discuss fiberings on a positive definite
$qc$-manifold $X$ with $qck$-group $\RR^3$. 
Section \ref{sub:Ca} reviews the the standard $qc$-structure $\sfD_0$ on 
the quaternionic Heisenberig Lie group $\cM$.
In Section \ref{producthyp} we construct a complete hyperK\"ahler metric
on the product $\RR\times X$ for a complete simply connected 
noncompact $qc$-manifold $X$. 
In Section \ref{abelnormal},
we determine the structure of $\Psh_{qc}(X)$ when $X$ has the $qck$-group
$\RR^3$.
It may be useful to determine $\Psh_{qc}(X)$
for several $qc$-manifolds like homogeneous spaces.
In Section \ref{Example1}, we classify spherical homogeneous $qc$-manifolds.
In Section \ref{app}, we calculate Ricci curvatures on 
$qc$-manifolds $X$ with $qck$-group $\RR^3$
using the associated Riemannian metric.
\vskip0.1cm

\section{Quaternionic contact group}\label{QContact group}
This section is mainly concerned with the equivalence of quaternionic contact structures.
Let $\HH$ be the field of quaternions $\{1,i,j,k\}$.
A \emph{quaternionic contact structure} 
is a codimension $3$-subbundle $\sfD$ on 
a $4n+3$-dimensional smooth manifold $X$ such that
$\displaystyle \sfD$ with 
$[\sfD,\sfD]$ generates $TX$. Moreover the following conditions are required:
There exists a non-degenerate ${\rm Im}\, \HH$-valued $1$-form 
$\om=\om_1i+\om_2j+\om_3k$ which represents $\sfD$, that is
$\displaystyle \ker\,\om=\mathop{\cap}_{\al=1}^3\ker\,\om_\al=\sfD$ 
such that
$\displaystyle \om\mathop{\we}\om\mathop{\we}\om\mathop{\we}\,  
\overbrace{d\om\we\cdots\we d\om}^n
\neq 0$ on $X$.
The $1$-form $\om$ is said to be a \emph{quaternionic contact form} on $X$.
Then the bundle of
endomorphisms $\{J_1,J_2,J_3\}$ 
defined by  
\begin{equation}\label{hypercom}
J_\ga=(d\om_\be|{\sf D})^{-1}\circ
(d\om_{\al}|{\sf D}):{\sf D} \ra {\sf D}\ \,((\alpha, \beta,\gamma)\sim(1,2,3))
\end{equation}
constitutes a hypercomplex structure on $\sfD$,\, 
that is $J_\al^2=-1$, $J_\al J_\be=J_\ga$.  
The Levi form 
$d\om_\al\circ J_\al: \sfD\times \sfD\lra\, \RR$
is a positive definite symmetric bilinear form on $\sfD$ (\cf Introduction).
Then $(X,\sfD,\om,\{J_k\}_{k=1}^{3})$ is called a 
 \emph{positive definite quaternionic contact} manifold
($qc$-manifold for short). See \cite{Biq}, \cite{IMV},
\cite{AK} for the definition and the reference therein.
Choose an ${\rm Im}\,\HH$-valued $1$-form $\om$ such that $\sfD=\ker\,\om$.
There is no canonical choice of $\om$ representing $\sfD$.
If $\om'$ represents $\sfD$, then it is easy to see that
there is a map ${v}\cdot b:X\ra \HH^*=\RR^+\times \Sp(1)$
such that $\om'=v\cdot b\,\om\bar b$. 

\begin{definition}\label{qcequivalence} 
Two ${\rm Im}\,\HH$-valued $1$-forms $\om$, $\om'$ are
\emph{qc-conformal} if $\om'=v\cdot b\, \om \bar b$
for some map $\displaystyle {v}\cdot b\in C^\infty(X,\RR^+\times \Sp(1))$.
\end{definition}
A quaternionic contact transformation is a diffeomorphism
$\al:X\ra\, X$ preserving $\sfD$ and the hypercomplex structure $\{J_\al\}_{\al=1}^3$.
More precisely, choose an ${\rm Im}\,\HH$-valued $1$-form $\om$ such that $\sfD=\ker\,\om$.
 The \emph{quaternionic contact group} 
$\Aut_{qc}(X)=\Aut_{qc}(X,\sfD,\, \{J_k\}_{k=1}^3\,)$ is defined by

\begin{equation}\label{Autqc}
\begin{split}\Aut_{qc}(X)=\bigl\{
\al\in\Diff(X)& \ \big|\  \al^*\om
=\lam_\al\cdot \om\cdot \overline{\lam_\al}
=u_\al\cdot a_\al\cdot \om\cdot \overline{a_\al}\,,\\
&\al_* \circ J_k=\sum_{j=1}^{3}
a_{kj}J_j\circ \al_*\ \, \mbox{on}\ \sfD\ \bigr\},
\end{split}\end{equation}for some map
$\lam_\al=\sqrt{u_\al}\cdot a_\al\in C^\infty(X,\RR^+\times \Sp(1))$ and 
the matrix $(a_{kj})\in C^\infty(X,\SO(3))$
is given by the conjugate of $a_\al$ on ${\rm Im}\,\HH$,
that is $\displaystyle z\to {a_\al}\cdot z\cdot \overline{a_\al}$.\\

This definition may depend on the choice of 
$\om$. 
Suppose $\om'=v\cdot b\cdot \om\cdot \bar b$\,
$({}^\exists\,\displaystyle {v}\cdot b: X\ra \RR^+\times \Sp(1))$.
The matrix $(b_{ij}):X\ra {\rm SO}(3)$ is defined by the conjugation 
$\displaystyle z\mapsto b\cdot z\cdot\bar b$.
Then the quaternionic structure $\{J'_1,J'_2,J'_3\}$
is obtained as 
\begin{equation}\label{eq.1}
d\om'_\be(J'_\ga \mbi{u},\mbi{v})=
d\om'_\al(\mbi{u},\mbi{v})\ \ \, ({}^\forall\, \mbi{u},\mbi{v}\in \sfD).
\end{equation}

\begin{lemma}\label{equiv}
Let $\displaystyle B=\left[\begin{array}{ccc}
b_{11} &b_{12}& b_{13}\\
b_{21} &b_{22}& b_{23}\\
b_{31} &b_{32}& b_{33}\\
\end{array}
\right]\in C^\infty(X,{\rm SO}(3))$ defined by the conjugate of $b\in C^\infty(X,\Sp(1)$.
Then 
\begin{equation}\label{conjuB}
\left[\begin{array}{c}
J_1'\\ J_2'\\ J_3'\\
\end{array}
\right]={}^tB\left[\begin{array}{c}
J_1\\ J_2\\ J_3\\
\end{array}
\right]\, .
\end{equation}
\end{lemma}

\begin{proof}
Since $d\om'_\be=v\sum_{j}b_{j\be}d\om_j$ on $\sfD$ as above,
using \eqref{eq.1}, it follows
$v\sum_{j}b_{j\be}d\om_j(J'_{\ga} \mbi{u},\mbi{v})=
v\sum_{j}b_{j\al}d\om_j(\mbi{u},\mbi{v})$. Thus
\begin{equation}\label{eq.2}\begin{split}
&\sum_{j}d\om_j((b_{j\be }J'_{\ga} -
b_{j \al })\mbi{u},\mbi{v})=0.\\
\end{split}
\end{equation}
Noting $d\om_1\circ J_1=d\om_2\circ J_2=d\om_3\circ J_3$,
\eqref{eq.2} is described: 
\begin{equation}\label{eq.3}\begin{split}
&d\om_1\circ J_1\bigl(-J_1(b_{1\be }J'_{\ga} -
b_{1\al })\mbi{u},\mbi{v})+
d\om_2\circ J_2(-J_2(b_{2\be }J'_{\ga} -
b_{2\al})\mbi{u},\mbi{v})\\
&\ \ \ \ +d\om_3\circ J_3(-J_3(b_{3\be }J'_{\ga} -
b_{3\al })\mbi{u},\mbi{v}\bigr)=0,\, \mbox{that is}\\
&d\om_1\circ J_1\bigl(J_1(b_{1\be }J'_{\ga}-b_{1\al })\mbi{u}+
J_2(b_{2\be }J'_{\ga}-b_{2\al })\mbi{u}\\
&\ \ \  \ \ \ \ \ \ \ \ \ +J_3(b_{3\be }J'_{\ga}-b_{3\al })\mbi{u},\mbi{v}\bigr)=0.\\
\end{split}
\end{equation}
By the non-degeneracy of $d\om_1\circ J_1$ for any $\mbi{v}\in\sfD$, we obtain
\[\bigl(J_1(b_{1\be }J'_{\ga}-b_{1\al })+
J_2(b_{2\be }J'_{\ga}-b_{2\al })+
J_3(b_{3\be }J'_{\ga}-b_{3\al })\bigr)\mbi{u}=0.\]
Equivalently
\begin{equation}\label{eq.4}\begin{split}
(J_1b_{1\be }+J_2 b_{2\be }+J_3b_{3\be })J'_{\ga}=
J_1b_{1\al }+J_2b_{2\al }+J_3b_{3\al}.
\end{split}\end{equation}Multiply \eqref{eq.4} by
$J_1b_{1\be },\,J_2b_{2\be},\, J_3b_{3\be }$ respectively:
\begin{equation*}\begin{split}
\bigl(-b_{1\be }^2+J_3b_{1\be }b_{2\be }-J_2b_{1\be }b_{3\be }\bigr)J'_{\ga}
&=-b_{1\be }b_{1\al }+b_{1\be }b_{2\al}J_3-a_{1\be }b_{3\al }J_2\\
\bigl(-J_3b_{1\be }b_{2\be }-b_{2\be }^2 +J_1b_{2\be }b_{3\be }\bigr)J'_{\ga}
&=-b_{2\be }b_{1\al }J_3-b_{2\be }b_{2\al }+b_{2\be }b_{3\al }J_1\\
\bigl(J_2b_{1\be }b_{3\be }-J_1b_{2\be }b_{3\be }-b_{3\be }^2\bigr)J'_{\ga}
&=b_{3\be }b_{1\al }J_2-b_{3\be }b_{2\al }J_1-b_{3\be }b_{3\al }\\
\end{split}\end{equation*}
Sum up these equations.
\begin{equation*}\begin{split}
-(b_{1\be }^2+b_{2\be }^2+b_{3\be }^2)
J'_r&=-(b_{1\be }b_{1\al }+b_{2\be }b_{2\al }+b_{3\be }b_{3\al })
+(b_{2\be }b_{3\al }-b_{3\be }b_{2\al })J_1\\
&\ \ \ +
(b_{3\be }b_{1\al }-b_{1\be }b_{3\al })J_2+
(b_{1\be }b_{2\al }-b_{2\be }b_{1\al })J_3.
\end{split}\end{equation*}Noting $\al\neq \be$,
$b_{1\be }^2+b_{2\be }^2+b_{3\be }^2=1,\,
b_{1\be }b_{1\al}+b_{2\be }b_{2\al }+b_{3\be }b_{3\al }=0$,
it follows
\begin{equation}\label{eq.5}\begin{split}
J'_r&=(b_{3\be }b_{2\al }-b_{2\be 2}b_{3\al })J_1
+(b_{1\be }b_{3\al }-b_{3\be }b_{1\al })J_2\\
&\ \ \ \ +
(b_{2\be }b_{1\al }-b_{1\be }b_{2\al })J_3.
\end{split}\end{equation}
As $(\al,\be,\ga)\sim (1,2,3)$,
we obtain the following equation from \eqref{eq.5}:
\begin{equation*}
\begin{split}
\left[\begin{array}{c}
J_1'\\
J_2'\\
J_3'
\end{array}
\right] =\left[\begin{array}{ccc}
b_{33}b_{22}-b_{23}b_{32} &b_{13}b_{32}-b_{33}b_{12}& b_{23}b_{12}-b_{13}b_{22}\\
b_{31}b_{23}-b_{21}b_{33} &b_{11}b_{33}-b_{31}b_{13}& b_{21}b_{13}-b_{11}b_{23}\\
b_{32}b_{21}-b_{22}b_{31} &b_{12}b_{31}-b_{32}b_{11}& b_{22}b_{11}-b_{12}b_{21}\\
\end{array}
\right]\left[\begin{array}{c}
J_1\\
J_2\\
J_3
\end{array}
\right]. 
\end{split}\end{equation*}
Represent ${}^tB=[\mbi{b}'_1,\mbi{b}'_2,\mbi{b}'_3]$ as 
column vectors. The above equation turns to
\begin{equation*}\begin{split}
\left[\begin{array}{c}
J_1'\\
J_2'\\
J_3'
\end{array}
\right]&=[\mbi{b}'_2\times \mbi{b}'_3,\mbi{b}'_3\times \mbi{b}'_1,
\mbi{b}'_1\times \mbi{b}'_2]
\left[\begin{array}{c}
J_1\\
J_2\\
J_3
\end{array}
\right] \ \mbox{and so},\\
\end{split}\end{equation*}
$\displaystyle \left[\begin{array}{c}
J_1'\\
J_2'\\
J_3'
\end{array}
\right]=[\mbi{b}'_1,\mbi{b}'_2,
\mbi{b}'_3]
\left[\begin{array}{c}
J_1\\
J_2\\
J_3
\end{array}
\right]={}^tB \left[\begin{array}{c}
J_1\\
J_2\\
J_3
\end{array}
\right]$,
the equation \eqref{conjuB}
follows.
\end{proof}

\begin{pro}\label{choice}
$\Aut_{qc}(X,\sfD,\om,\{J_k\}_{k=1}^3)=
\Aut_{qc}(X,\sfD,\om',\{J'_k\}_{k=1}^3)$ 
whenever $\om'=v\cdot b\cdot \om\cdot \bar b$, that is, it is uniquely determined by
the qc-conformal class of $\om$.
Henceforth the definition of $\displaystyle \Aut_{qc}(X,\sfD)$ makes sense.
\end{pro}

\begin{proof}
Let $\al\in \Aut_{qc}(X,\sfD,\om,\{J_k\}_{k=1}^3)$.
By the definition of \eqref{Autqc} for $\om$,

\begin{enumerate}
\item[(1)] $\displaystyle \al^*\om=u_\al\cdot a_\al\cdot \om\cdot \overline{a_\al}$.
\item[(2)]
 $\displaystyle \al_* \left[
J_1\, J_2\, J_3
\right]
=\overline{a_\al}
\left[
J_1\, J_2\, J_3
\right]a_\al\cdot \al_*$.
\end{enumerate}
From Lemma \ref{equiv} note
$\displaystyle\left[
J'_1\,J'_2\,J'_3
\right]
=b\left[
J_1\,J_2\,J_3
\right]
\bar b$.
Suppose $\om'=v\cdot b\cdot \om\cdot \bar b$.
Then $\displaystyle\al^*\om'=\al^*v\cdot \al^*b\cdot \al^*\om\cdot \overline{\al^*b}
=(\al^*v\cdot u_\al\cdot v^{-1})\cdot
(\al^*b\cdot a_\al\cdot \bar b)\cdot \om'\cdot
(\overline{\al^*b\cdot a_\al\cdot \bar b})$. 
\begin{equation*}\begin{split}
\al_* \left[
J'_1\,J'_2\,J'_3
\right]
&=\al_*\left(b
\left[
J_1\,J_2\,J_3\right]
\bar b\right)
=b\al_*\left[J_1\,J_2\,J_3
\right]\bar b\\
&=
b\cdot\overline{a_\al}\left[
J_1\,J_2\,J_3
\right]a_\al\cdot \bar b \cdot \al_*\\
&=(b\cdot\overline{a_\al}\cdot\overline{\al^*b})\cdot\left[
J'_1\,J'_2\,J'_3 \right]\cdot (\al^*b\cdot \al_\al\cdot \bar b)\\
&=(b\cdot\overline{a_\al}\cdot\overline{\al^*b})\cdot\left[
J'_1\,J'_2\,J'_3
\right]\cdot (\overline{b\cdot\overline{a_\al}\cdot\overline{\al^*b}}). 
\end{split}\end{equation*}
Noting $\displaystyle 
\overline{b\cdot\overline{a_\al}\cdot\overline{\al^*b}\,}=
{\al^*b\cdot a_\al\cdot \bar b}$,
this implies $\al\in \Aut_{qc}(X,\sfD,\om',\{J'_k\}_{k=1}^3)$ by the definition.
Conversely let $\al\in \Aut_{qc}(X,\sfD,\om',\{J'_k\}_{k=1}^3)$.
\begin{enumerate}
\item[(1')] $\displaystyle \al^*\om'=u'_\al\cdot a'_\al\cdot \om'
\cdot \overline{a'_\al}$.
\item[(2')]
 $\displaystyle \al_* \left[
J'_1\,J'_2\,J'_3\right]
=\overline{a'_{\al}}
\left[
J'_1\,J'_2\,J'_3\right]
{a'_\al}\cdot \al_*$.
\end{enumerate}
Similarly as above, 
\begin{equation*}\begin{split}
\al^*\om=\al^*(v^{-1}\cdot \bar b\cdot \om'\cdot b)
=(\al^*v^{-1}\cdot u'_\al\cdot v)\cdot (\overline{\al^*b}\cdot a'_\al\cdot b)
\cdot \om\cdot (\overline{\overline{\al^*b}\cdot a'_\al\cdot b}),\\
\end{split}\end{equation*}
\begin{equation*}\begin{split}
\al_* \left[
J_1\,J_2\,J_3\right]
&=\al_*\left(\bar b
\left[J'_1\,J'_2\,J'_3\right]b\right)=
\bar b\cdot
\al_*\left[
J'_1\,J'_2\,J'_3
\right]\cdot b\\
&=
(\bar b\cdot\overline{a'_\al}\cdot \al^*b)
\cdot\left[J_1\, J_2\, J_3\right]
\cdot (\overline{\al^*b}\cdot{a'_\al}\cdot b)\cdot \al_*\\
&=\overline{(\overline{\al^*b}\cdot
a'_\al\cdot b)}\cdot\left[J_1\, J_2\, J_3\right]
(\overline{\al^*b}\cdot a'_{\al}\cdot b)\cdot \al_*.
\end{split}\end{equation*}
Hence
$\displaystyle \al\in \Aut_{qc}(X,\sfD,\om,\{J_k\}_{k=1}^3)$. It follows
$$ \Aut_{qc}(X,\sfD,\om,\{J_\al\}_{\al=1}^3)=\Aut_{qc}(X,\sfD,\om',\{J'_\al\}_{\al=1}^3).$$
\end{proof} 

\section{Conformal invariant}\label{Confinv}
Let $\displaystyle \Aut_{qc}(X)=\Aut_{qc}(X,\sfD,\,\{J_\al\}_{\al=1}^3\,)$.
If $C^\infty(X,\RR^+)$ is the module
consisting of smooth functions of $X$ to the positive numbers $\RR^+$,
then it is
endowed with an action of $\al\in \Aut_{qc}(X)$.
For $f\in C^\infty(X,\RR^+)$, 
\[
(\al_* f)(x)=f(\al^{-1}x)\,\ (x\in X).\]
Thus $C^\infty(X,\RR^+)$ is a smooth $G$-module (Compare \cite{OK}.)
For a form $\om=\om_1i+\om_2j+\om_3k$, denote the norm as
$\displaystyle |\om|=\sqrt {\om_1^2+\om_2^2+\om_3^2}=\sqrt{\langle \om,\om\rangle}$.

\begin{theorem}[\cite{OK}]\label{vanish}
Let $G$ be a closed subgroup of $\Aut_{qc}(X)$.
There is a cocycle $\mu_{qc}\in H^1_d(G,C^\infty(X,\RR^+))$
which is a \emph{$qc$-conformal invariant}.
\end{theorem}

\begin{proof}
Let $\al\in \Aut_{qc}(X)$ such that $\displaystyle
\al^*\om=u_\al\cdot a_\al \om \cdot\overline{a_\al}$.
For $\al, \be\in \Aut_{qc}(X)$, calculate
\begin{equation*}\begin{split}
(\al\be)^*\om&=u_{\al\be}\cdot a_{\al\be} \om \cdot\overline{a_{\al\be}},\\
\be^*\al^*\om&=\be^*(u_\al\cdot a_\al \om\cdot \overline{a_\al})
=(\be^*u_\al\cdot u_\be)\cdot (\be^*a_\al\cdot  a_\be)\om\cdot
(\overline{\be^*a_\al\cdot a_\be}).\\
\end{split}
\end{equation*}
(Note that $\displaystyle \be^*\overline{a_\al}(x)=
\overline{a_\al}(\be x)=\overline{a_\al(\be x)}=
\overline{\be^*a_\al}(x)$).
Since each $a_\ga$ $(\ga\in G)$ is a map from $X$
to $\Sp(1)$, we have
$\displaystyle |(\al\be)^*\om|=u_{\al\be}|\om|=\be^*u_\al\cdot u_\be|\om|$
by the above equation.
Thus the smooth maps $u_\al, u_\be, u_{\al\be} : X\ra \RR^+$ satisfy
\begin{equation}\label{crossed}
 u_{\al\be}=\be^*u_\al\cdot u_\be\ \ \mbox{on}\ X. 
\end{equation}
As in \cite{OK}, 
define $\lam_{\om}: G\ra C^\infty(X,\RR^+)$ to be
\begin{equation}\label{crossedhom}
\lam_{\om}(\al)=\al_*u_\al.
\end{equation}
In particular, $\displaystyle \lam_{\om}(\al)(x)=u_\al(\al^{-1}x)$.
It suffices to show $\lam_{\om}$ is a crossed homomorphism.
Calculate
\begin{equation*}
\begin{split}
\lam_{\om}(\al\be)(x)&=(\al\be)_*u_{\al\be}(x)=u_{\al\be}(\be^{-1}\al^{-1}x)\\
 &=\be^*u_\al(\be^{-1}\al^{-1}x)\cdot u_\be(\be^{-1}\al^{-1}x)\ \, (\eqref{crossed})\\
&= \lam_{\om}(\al)(x)\cdot\al_*\lam_{\om}(\be)(x)
=(\lam_{\om}(\al)\cdot\al_*\lam_{\om}(\be))(x).
\end{split}\end{equation*}Hence
$\displaystyle \lam_{\om}(\al\be)=\lam_{\om}(\al)\cdot\al_*\lam_{\om}(\be)$,
that is $\lam_{\om}$ is a crossed homomorphism so that
$\displaystyle [\lam_{\om}]\in H^1_d(G,C^\infty(X,\RR^+))$.
We show $[\lam_{\om}]$ is a quaternionic conformal invariant.
Suppose $\om'=u\cdot b\,\om \bar b$ for some positive function $u$ on $X$ and
a map $b:X\ra \Sp(1)$.
For $\al\in G$, let 
$\displaystyle \al^*\om'=u'_\al\cdot a'_{\al}\om'\, \overline{a'_{\al}}$.
Then $\displaystyle \al^*\om'=u'_\al\cdot a'_{\al}u\cdot b\,\om \bar b\, \overline{a'_{\al}}
=(u'_\al\cdot u)\cdot (a'_{\al}\cdot b)\om\cdot (\overline{a'_{\al}\cdot b})$,
while
\begin{equation*}\begin{split}
\al^*\om'&=\al^*(u\cdot b\,\om \bar b)
=(\al^*u\cdot u_\al)\cdot (\al^*b\cdot a_\al)\om\cdot
(\overline{\al^*b\cdot a_\al}).
\end{split}\end{equation*}
Taking the norm $|\al^*\om'|$, it follows
$\displaystyle u'_\al\cdot u=\al^*u\cdot u_\al$,\ that is\,
$\displaystyle u'_\al=\al^*u\cdot u^{-1}\cdot u_\al$
on $X$.
For $\al^{-1}x\in X$, as $\displaystyle u'_\al(\al^{-1}x)=\al^*u(\al^{-1}x)\cdot u^{-1}(\al^{-1}x)\cdot u_\al(\al^{-1}x)$, it follows $\displaystyle 
\al_*u'_\al(x)=u(x)\cdot (\al_*u)^{-1}(x) \al_*u_\al(x)$, that is
$\displaystyle 
\al_*u'_\al(x)\cdot (\al_*u)(x)\cdot u(x)^{-1}= \al_*u_\al(x)$ $({}^\forall\, x\in X)$. This shows
$\displaystyle \al_*u'_\al\cdot \delta^0(u)(\al)=\al_*u_\al$. 
Hence
$\displaystyle [\lam_{\om}]=[\lam_{\om'}]$ and so
$[\lam_{\om}]$ is a quaternionic conformal invariant.
We may put 
\begin{equation}\label{muinv}
\mu_{qc}=[\lam_{\om}]\in H^1_d(G,C^\infty(X,\RR^+)).\end{equation}
\end{proof}

\begin{remark}
$(i)$\, It is shown in \cite{OK} that
$\displaystyle H^i_d(G,C^\infty(X,\RR^+))=0$ 
$(i\geq 1)$
provided that $G$ acts properly on $X$. In particular, $\mu_{qc}=0$.
\noindent
$(ii)$\, Since we evaluate the norm of $\om$,
$\lam_{\om}$ is in fact a \emph{conformal invariant}
whenever $\om'=u\cdot \om$.
\end{remark}

\begin{definition}\label{almosth}
The \emph{$qc$-Hermitian group} is denoted by 
\begin{equation}\label{qcPsh}
\begin{split}\Psh_{qc}(X, \om,\{J_k\}_{k=1}^3)=\bigl\{
\al\in\Diff(X)
\, \big|\  \al^*\om
&=a_\al \cdot\om\cdot \overline{a_\al}\,,\\
\al_* \circ J_k&=\sum_{j=1}^{3}
a_{kj}J_j\circ \al_*|_{\sfD}\ \},\
\end{split}\end{equation}where
$a_\al\in C^\infty(X, \Sp(1))$ induces
$(a_{ij})\in C^\infty(X,{\rm SO}(3))$ as
its conjugate. 
\end{definition}\noindent(Compare \ref{Autqc}.)
By Definition \ref{Autqc} note 
\begin{equation*}
 \Psh_{qc}(X, \om,\{J_k\}_{k=1}^3)
\leq \Aut_{qc}(X,\sfD,\{J_k\}_{k=1}^3).
\end{equation*}
Using the equality 
$\displaystyle d\om_1\circ J_1=
d\om_2\circ J_2=d\om_3\circ J_3$ on $\sfD$ (\cf \eqref{eq:IJK}, also \cite{AK}),
each $\al\in \Psh_{qc}(X,\om,\{J_k\}_{k=1}^3)$ satisfies
\begin{equation*}
\begin{split}
&3d\om_1(J_1 \al_*\mbi{a},\al_*\mbi{b})=
\sum_{j=1}^{3}d\om_j(J_j \al_*\mbi{a},\al_*\mbi{b})
=\sum_{j,k}d\om_j(\al_*a_{jk}J_k\mbi{a},\al_*\mbi{b})\\
&=\sum_{j,k}d\al^*\om_j(a_{jk}J_k\mbi{a},\mbi{b})=\sum_{j,k,\ell}d\om_\ell a_{\ell j}(a_{jk}J_k\mbi{a},\mbi{b})\\
&=\sum_{k,\ell}^{}d\om_\ell(\delta_{\ell k}J_k\mbi{a},\mbi{b})
=\sum_{k=1}^{3}d\om_k(J_k\mbi{a},\mbi{b})
=3d\om_1(J_1\mbi{a},\mbi{b}),\ \mbox{that is}
\end{split}\end{equation*}
\begin{equation}\label{3equiv}
 d\om_1(J_1 \al_*\mbi{a},\al_*\mbi{b})=
d\om_1( J_1 \mbi{a},\mbi{b})\, \ ({}^\forall\,\mbi{a},\mbi{b}\in \sfD).
\end{equation}

Since $d\om_1\circ J_1 : \sfD\times \sfD\ra \RR$
defined by $d\om_1(J_1 \mbi{a}, \mbi{b})$ 
is positive definite in our case, 
a $qc$-manifold $(X, \sfD\,(=\ker\, \om),\{J_\al\}_{\al=1}^3)$
assigns to $\sfD$ a Riemannian metric (\cf \eqref{eq:IJK}) 
\begin{equation}\label{Riemann}
g_\om(\mbi{x},\mbi{y})=\sum_{i=1}^{3}\om_i(\mbi{x})\cdot \om_i(\mbi{y})+
d\om_1(J_1 \mbi{x}, \mbi{y})\, \ \ ({}^\forall\, \mbi{x}, \mbi{y}\in TX).
\end{equation} 
The isometry group $\Isom(X,g_\om)$ is not related to 
$\Aut_{qc}(X)$ in general.

\begin{pro}\label{pqc}
$\displaystyle \Psh_{qc}(X,\om,\{J_\al\}_{\al=1}^3)\leq \Isom(X,g_\om)$.
Then $\Psh_{qc}(X,\om)$ acts properly on $X$.
When $X$ is compact, $\Psh_{qc}(X,\om,\{J_\al\}_{\al=1}^3)$ is a 
compact Lie group.
\end{pro}

\begin{proof}
If $\al\in \Psh_{qc}(X,\om,\{J_\al\}_{\al=1}^3)$, then
$\al^*\om=a_\al\cdot \om\cdot \overline{a_\al}$ by the definition.
As $\sum_{i=1}^{3}\om_i\cdot \om_i=-\om\cdot \om\in \RR$,
it follows $\displaystyle \al^*\om\cdot \al^*\om=a_\al\cdot \om \cdot\overline{a_\al}\cdot
a_\al\cdot \om \cdot \overline{a_\al}=\om\cdot\om$.
By \eqref{3equiv}, 
$\displaystyle \al^*(d\om_1\circ J_1)=d\om_1\circ J_1$.
This implies $\displaystyle \al^*g=g$ such that 
$\displaystyle \Psh_{qc}(X,\om,\{J_\al\}_{\al=1}^3)\leq \Isom(X,g_\om)$.
\end{proof}

\begin{theorem}\label{qvanish}
Let $G$ be a closed subgroup of 
$\Aut_{qc}(X,\sfD,\, \{J_k\}_{k=1}^3\,)$.
Then $\mu_{qc}=0$ in $H^1_d(G, C^\infty(X,\RR^+))$
if and only if $G$ acts properly on $X$.
In that case, there is an ${\rm Im}\,\HH$-valued $1$-form
$\eta$ \emph{conformal} to $\om$, such that  
\[G\leq \Psh_{qc}(X,\eta,\{J_k\}_{k=1}^3\,).\]
\end{theorem}

\begin{proof}
Suppose $\mu_{qc}=[\lam_{\om}]=0$ in $H^1_d(G,C^\infty(X,\RR^+))$
for an ${\rm Im}\,\HH$-valued $1$-form $\om$.
Then $\lam_{\om}=\delta^0v$ for $v\in C^\infty(X,\RR^+)$.
It follows $\displaystyle \al_*u_\al(x)=\al_*v(x)v(x)^{-1}$\ \, $({}^\forall\, \al\in G)$.
So $u_\al(\al^{-1}x)v(x)=v(\al^{-1}x)$. In particular,
$\displaystyle u_\al(x)v(\al x)=v(x)$ $({}^\forall\, x\in X)$,
that is $\displaystyle u_\al\cdot \al^*v=v$.
If $\al\in G$, then $\al^*\om=u_\al\cdot a_\al\cdot \om\,\cdot \overline{a_\al}$ as above.
Put 
\begin{equation}\label{conf-om}
 \eta=v\cdot\om.
\end{equation}
Then
$\displaystyle \al^*\eta=\al^*v\cdot u_\al\cdot a_\al\cdot \om\cdot\, \overline{a_\al}
=v\cdot a_\al\cdot \om\cdot\, \overline{a_\al}=a_\al\cdot \eta\cdot \overline{a_\al}$.
Moreover,
if $\{J'_j\}_{j=1}^3$ is a hypercomplex structure for $\displaystyle
\eta=\eta_1i+\eta_2j+\eta_3k$,
then $\displaystyle d\eta_i(J'_j \mbi{a},\mbi{b})=d\eta_k(\mbi{a},\mbi{b})$
$({}^\forall\, \mbi{a}, \mbi{b}\in \sfD)$ by the definition.
It implies $\displaystyle vd\om_i(J'_j \mbi{a},\mbi{b})=vd\om_k(\mbi{a},\mbi{b})$
on $\sfD$, thus $\displaystyle d\om_i(J'_j \mbi{a},\mbi{b})=d\om_k(\mbi{a},\mbi{b})
=d\om_i(J_j \mbi{a},\mbi{b})$. Replace $\mbi{b}$ by $J_i\mbi{b}$.
The non-degeneracy of $d\om_i\circ J_i$ on $\sfD$ $(i=1,2,3)$
implies $J'_j=J_j$ $(j=1,2,3)$. 
Hence
$\displaystyle \al\in\Psh_{qc}(X,\eta,\{J_\al\}_{\al=1}^3)$,
 that is $G\leq \Psh_{qc}(X,\eta,\{J_\al\}_{\al=1}^3)$.
By Proposition \ref{pqc}, $\displaystyle \Psh_{qc}(X,\eta,\{J_\al\}_{\al=1}^3)$
acts properly on $X$, so does $G$ .

Conversely, if $G$ acts properly on $X$, then it follows from \cite[Theorem 14]{OK1} that
$H^1_d(G,C^\infty(X,\RR^+))=0$.
In particular, $\mu_{qc}=0$ by Theorem \ref{vanish}.
\end{proof}

\section{ $3$-dimensional complementary distribution $V$}\label{complemintg}
When $(X, \om,\{J_\al\}_{\al=1}^3)$ is a $qc$-manifold,
there is a $qc$-structure $(\eta,\{J_\al\}_{\al=1}^3)$ $qc$-conformal to $\om$ by Theorem A. Let $\eta=\eta_1i+\eta_2j+\eta_3k$ be the $1$-form
with $\sfD= \ker\, \eta$.
As $\displaystyle d\eta_\al:\sfD\times \sfD\,\ra\, \RR$ is non-degenerate,
there exist three vector fields
$\xi_\al$ (\cf \cite{Biq}) such that
\begin{equation*}
\eta_\al(\xi_\be)=\delta_{\al\be},\,\
d\eta_\al(\xi_\al,\mbi{u})=0 \,
\ (\al,\be=1,2,3, \, {}^\forall\, \mbi{u}\in \sfD).\\
\end{equation*}

Put $\displaystyle V=\{\xi_1,\xi_2,\xi_3\}$ associated with $\eta$.
Then the above equations determine $V$ uniquely. 

\begin{definition}\label{qcV}
 Let $V$ be the $3$-dimensional distribution 
on a $qc$-manifold $(X,\eta,\{J_\al\}_{\al=1}^3)$
with $\displaystyle \sfD=\ker\, \eta$. 
$(${\rm \cf}$\eqref{fix})$.
If\ $V$ generates a subgroup $\cR$ of $\Psh_{qc}(X,\eta,\{J_\al\}_{\al=1}^3)$,
then $V$ is said to a
$qck$-distribution on $X$ associated with $\eta$. $\cR$ is called a $qck$-group.
\end{definition}
Suppose $V$ is a $qck$-distribution on $X$. Noting 
$[\xi_\al,\mbi{u}]\in \sfD$ $(\al=1,2,3)$, the above equation shows  
\begin{equation}\begin{split}\label{ijk}
d\eta_\al(\xi_\be,\mbi{u})=0\
\ (\al,\be=1,2,3, \, {}^\forall\, \mbi{u}\in \sfD).\\
\end{split}
\end{equation}

There is a decomposition:
\begin{equation}\label{fix}
TX=V\oplus\sfD.
\end{equation}

Denote $J_\al\xi_\be=0$ $(\al,\be=1,2,3)$ as usual.
Let $\mbi{x}=\sum_\al a_{\al\be}\xi_\al+\mbi{u},\,
\mbi{y}=\sum_\al b_{\al\be}\xi_\al+\mbi{v}$ as above.
For any such $\mbi{x},\mbi{y}\in TX$, it follows
$\displaystyle d\eta_\al(J_\al \mbi{x},\mbi{y})=d\eta_\al(J_\al\mbi{u},\mbi{v})$
by \eqref{ijk}.
In particular, from \eqref{hypercom}
\begin{equation}\label{eq:IJK}
d\eta_1(J_1 \mbi{x},\mbi{y})=d\eta_2(J_2 \mbi{x},\mbi{y})=d\eta_3(J_3 \mbi{x},\mbi{y})\ \ ({}^\forall\, \mbi{x}, \mbi{y}\in TX).
\end{equation}

\begin{pro}\label{isoqc}
Let $V=\{\xi_1,\xi_2,\xi_3\}$ be the
$3$-dimensional $qck$-distribution on $X$.
If $h\in \Psh_{qc}(X, \sfD, \eta,\, \{J_\al\}_{\al=1}^3)$, then
\vskip0.2cm
\begin{enumerate}
\item[(i)]
$\displaystyle d\eta_1(J_1 h_*\mbi{x}, h_*\mbi{y})=
d\eta_1(J_1 \mbi{x}, \mbi{y})$\ \, $({}^\forall\, \mbi{x},\mbi{y}\in TX)$.
\item[(ii)]
$h_*V=V$. 
\end{enumerate}
\end{pro}
\begin{proof}
(1) First similarly to the equation \eqref{3equiv}, we obtain 
$$ d\eta_1(J_1 h_*\mbi{u}, h_*\mbi{y})=
d\eta_1(J_1\mbi{u}, \mbi{y})\ \, ({}^\forall\, \mbi{u}\in \sfD, {}^\forall\,
\mbi{y}\in TX).$$

\noindent Then we prove (ii) $h_*V=V$.
By (1),
$\displaystyle d\eta_1(J_1 h_*\mbi{u},h_*\xi_\al)=
d\eta_1( J_1 \mbi{u},\xi_\al)=0$ for ${}^\forall\,\mbi{u}\in \sfD$.
 Suppose $h_*\xi_\al=\sum b_{\al\be}\xi_\be+\mbi{v}$
for some $\mbi{v}\in \sfD$. 
Since $\displaystyle d\eta_1(J_1 h_*\mbi{u},h_*\xi_\al)=
d\eta_1(J_1 h_*\mbi{u},\sum b_{\al\be}\xi_\be+\mbi{v})
=d\eta_1(J_1 h_*\mbi{u},\mbi{v})$,
the non-degeneracy
of $d\eta_1\circ J_1$ on $\sfD$ implies $\mbi{v}=\mbi{0}$.
Hence $h_*\xi_\al\in V$.

We prove (i) for any $\mbi{x},\mbi{y}\in TX$. 
If $\displaystyle \mbi{x}=\sum a_{\al\be}\xi_\al+\mbi{u}$ as above, then
$\displaystyle h_*\mbi{x}=\sum a_{\al\be}h_*\xi_\al+h_*\mbi{u}$ where 
$h_*\xi_\al\in V$ by (ii). Thus
$\displaystyle d\eta_1(J_1 h_*\mbi{x},h_*\mbi{y})=d\eta_1(J_1 h_*\mbi{u},h_*\mbi{y})
=d\eta_1(J_1 \mbi{u},\mbi{y})=d\eta_1(J_1 \mbi{x},\mbi{y})$.
\end{proof}

In particular, there is an analogy of
the Killing (Reeb) fields of Sasaki $CR$-manifolds to $qc$-manifolds from Proposition \ref{isoqc}.
We may put $\displaystyle \Psh_{qc}(X)=\Psh_{qc}(X, \sfD, \eta,\{J_\al\}_{\al=1}^3)$
for short. By Proposition \ref{isoqc} it follows
\begin{corollary}\label{normalizer}
Suppose a $qck$-distribution $V$ generates a subgroup $\cR\leq \Psh_{qc}(X)$.
If $\displaystyle N_{\Psh_{qc}(X)}(\cR)$
is the normalizer of $\cR$ in $\Psh_{qc}(X)$, then 
$$\Psh_{qc}(X)=N_{\Psh_{qc}(X)}(\cR).$$
\end{corollary}

\begin{remark}
Even though $\om$ generates the $qck$-distribution,
any $\eta$ $qc$-conformal to $\om$ does not necessarily generate a $qck$-distribution.
\end{remark}

Suppose a smooth $4n+3$-dimensional \emph{compact} manifold $M$ admits
a $qc$-structure $(\sfD,\{J_\al\}_{\al=1}^3)$.
By Theorem A, there exists 
an ${\rm Im}\,\HH$-valued $1$-form $\eta$ such that
$\ker\, \eta=\sfD$. Then    
$\Aut_{qc}(M,\sfD,\{J_\al\}_{\al=1}^3)=\Psh_{qc}(M,\eta,\{J_\al\}_{\al=1}^3)$
unless $M$ is $qc$-conformal to $S^{4n+3}$.
Suppose $V$ is a $qck$-distribution on $M$.
$V$ generates a $qc$-subgroup $T\leq
\Psh_{qc}(M)$ (\cf (2) Definition \ref{qcV}).
As $\Psh_{qc}(M)$ is a \emph{compact} Lie group by Proposition \ref{pqc},
so is the closure $\displaystyle \bar {T}$ of $T$.
Since $\displaystyle \bar {T}$ is connected,
it follows from \eqref{fix} that
$\displaystyle {\bar T}=T$.
Then $T$ is isomorphic to either $\Sp(1)$
or $T^3$. The compact case falls into the following dichotomy.

\begin{pro}\label{class3} Let $(M,\eta,\{J_\al\}_{\al=1}^3)$ be a \emph{compact}
$qc$-manifold with the $qck$-distribution $V=\{\xi_1,\xi_2,\xi_3\}$.
If $T\leq\Psh_{qc}(M)$ is
the subgroup generated
by $V$, then either one of the following
occurs exactly:
\begin{enumerate}
\item $T$ is isomorphic to $\Sp(1)\leq \Psh_{qc}(M)$.
$M$ is a quaternionic $CR$-manifold $(3$-Sasaki manifold\,$)$ for which
$\displaystyle V$
coincides with the common kernel of $\rho_\al=d\eta_\al+2\eta_\be\wedge\eta_\ga$\,
$(\al=1,2,3)$, $(\al,\be,\ga)\sim (1,2,3)$;
$$\ \ \ \ \, V=\{\xi\mid \rho_\al(\xi,\mbi{x})=0,\,
{}^\forall\, \mbi{x}\in TM\ (\al=1,2,3)\,\}.$$
\item $T$ is isomorphic to $T^3\leq \Psh_{qc}(M)$.
$M$ is a $T^3$-orbibundle over a hyperK\"ahler orbifold $Z$.
In this case,  \[ V=\{\xi\mid d\eta_\al(\xi,\mbi{x})=0, \,
{}^\forall\, \mbi{x}\in TM\, (\al=1,2,3)\,\}.\]
\end{enumerate}
\end{pro}

\begin{proof}
(1)\ If $T=\Sp(1)$, then $\displaystyle V=\{\xi_\al,\xi_\be,\xi_\ga\}$
$((\al,\be,\ga)\sim (1,2,3))$
satisfies $[\xi_\al,\xi_\be]=2\xi_\ga$.
If $\ker\,\rho_\al=\{\xi\in TM \mid
\rho_\al(\xi,\mbi{x})=0\ {}^\forall\, \mbi{x}\in TM\}$,
then $V=\mathop{\cap}_{\al=1}^3\ker \rho_\al$. $M$ is
a quaternionic $CR$-manifold (\cf \cite{AK}). 
$(2)$\  
If $T=T^3$, then
$V=\{\xi_1,\xi_2,\xi_3\}$ is the commutative algebra and  so
$d\om_\ga(\xi_\al,\mbi{x})=0$ $({}^\forall\, \mbi{x}\in TM)$.
If $p:M\ra\, Z=M/T^3$ is the projection, then this equation gives a hyperK\"ahler form $\Om_\al$ on the orbifold $Z$ such that $p^*\Om_\al=d\om_\al$, 
(see \eqref{hKaehler} for the detail, also \cite{KA2}).
\end{proof}
\medskip

\begin{remark}
If $M$ is a $qCR$ $(3$-Sasaki$)$ manifold, then
for each $S^1\leq S^3$
the orbifold $M/S^1$ is the twistor space over the quaternionic K\"ahler orbifold
$\displaystyle S^2\ra\, M/S^1\lra\,  M/S^3$. Moreover,
it is shown in \cite{BG1}, \cite{IK}
 that $M/S^1$ is a
K\"ahler-Einstein orbifold of positive scalar curvature.
For $(2)$, $M$ is a $qc$-Einstein manifold with zero $qc$-scalar curvature
by the result of \cite{IMV1}. 
\end{remark}

\section{Non-compact $qc$-manifold with $qck$-distribution $\RR^3$}\label{sec:fib}
Since a $qc$-manifold $X$ is non-compact simply connected in our case,
we may assume the $3$-dimensional $qck$-subgroup $\cR\leq \Psh_{qc}(X)$ 
is also a non-compact simply connected Lie group. 
From Definition \ref{almosth},
$\displaystyle r^*\eta=a_r\cdot \eta\cdot \overline {a_r}$ for $r\in \cR$
where $a_r\in C^\infty(X,\Sp(1))$ in general. We note
the following.

\begin{pro}\label{soR3}
If each $a_r$ is a constant map, that is $a_r$ is
an element of $\Sp(1)$, then $\cR=\RR^3$. In addition
$t^*\eta=\eta$,\, $t_*J_\al=J_\al t_*$ for any $t\in \RR^3$.
\end{pro}

\begin{proof}
The correspondence $\displaystyle
\nu(r)=a_r$ is a homomorphism of $\cR$ to 
$\Sp(1)$ ($[a_r]\in \SO(3)$ if necessary). If $\cR$ is semisimple, then it is isomorphic to
the universal covering $\displaystyle \widetilde{\SL(2,\RR)}$.
As $\Sp(1)$ is simply connected, $\nu:\cR\ra \Sp(1)$
turns to an isomorphism, which is impossible by our hypothesis,
or $\nu(\cR)=\{1\}$. If $\cR$ has the nontrivial radical,
then $\cR$ is solvable. It follows 
either $\nu(\cR)=\{1\}$ or $\nu(\cR)=S^1\leq \Sp(1)$. In case
$\nu(\cR)=S^1$,  
the identity component of $\ker\,\nu$ has $2$-dimension. 
Let $V=\{\xi_\al,\xi_\be,\xi_\ga\}$ be the $qck$-distribution which generates $\cR$.
We may assume $(\ker\,\nu)^0$ induces $\displaystyle
\{\xi_\be,\xi_\ga\}$ for example.
Noting $r^*\eta=\eta$ for $r\in (\ker\,\nu)^0$,
it follows $\displaystyle\cL_{\xi_\be}\eta=0$.
By \eqref{ijk}, $\displaystyle d\eta(\xi_{\be}, \mbi{x})=0$\,
 $({}^\forall\, \mbi{x}\in TX)$.
Taking $\mbi{x}=\xi_\al$, 
$\displaystyle \eta([\xi_{\be},\xi_\al])=0$.
Then $[\xi_{\be},\xi_\al]\in \sfD=\ker\,\eta$. 
Since $[\xi_k,J_1\mbi{v}]\in \sfD$ by \eqref{ijk}
$(k=\al,\be,\,{}^\forall\, \mbi{v}\in \sfD)$,
Jacobi identity implies
$\displaystyle d\eta_1([\xi_\be,\xi_\al],J_1\mbi{v})=0$.
The non-degeneracy of $d\eta_1\circ J_1$
shows $[\xi_\be,\xi_\al]=\mbi{0}$.
Similarly taking $\mbi{x}=\xi_\ga$, we have $[\xi_\be,\xi_\ga]=\mbi{0}$.
As $\cL_{\xi_\ga}\eta=0$ also, $[\xi_\ga,\xi_\al]=\mbi{0}$ as above.
Hence $V$ generates $\cR=\RR^3$ which shows $\cL_{\xi_\al}\eta=0$.
As a matter of fact, $\nu(\cR)=\{1\}$.
When $\cR$ is simply connected semisimple, $\nu(\cR)=\{1\}$ 
so the above argument applies to show $\cR=\RR^3$ which were
impossible. 
In addition, let $\displaystyle
d\eta_\be(J_\ga(t_*\mbi{u}),t_*\mbi{v})=d\eta_{\al}(t_*\mbi{u},t_*\mbi{v})$
$({}^\forall\, \mbi{u},\mbi{v}\in {\sfD})$ by \eqref{hypercom}.
Since $t^*\eta=\eta$,
$\displaystyle
d\eta_\be(t^{-1}_*J_\ga t_*\mbi{u},\mbi{v})=d\eta_{\al}(\mbi{u},\mbi{v})=
d\eta_\be(J_\ga\mbi{u},\mbi{v})$. The non-degeneracy of $d\eta_\be$
implies $t^{-1}_*J_\ga t_*\mbi{u}=J_\ga\mbi{u}$,
that is $J_\al t_*=t_*J_\al$ on $\sfD$ $(\al=1,2,3)$.
\end{proof}

\begin{remark}\label{posibty}
When $a_r\in C^\infty(X,\Sp(1))$ is a nontrivial smooth map,
we do not know whether
a solvable group $(\neq \RR^3)$ or
$\displaystyle \widetilde{\SL(2,\RR)}$
may occur as $\cR$. \end{remark}

We discuss the fiberings of a non-compact simply connected
 $qc$-manifold $X$ with a $qck$-group $\RR^3$.
Let $TX=V \oplus \sfD$ where $V=\langle \xi_1,\xi_2,\xi_3\rangle$
is the $qck$-distribution.
Take $J_1$ from $\displaystyle \{J_\al\}_{\al=1}^3$ on $\sfD$.
Put $\displaystyle E=\langle \xi_2,\xi_3\rangle \oplus \sfD$.
Define an almost complex structure $\bar J_1$ on the distribution 
$E$ to be
\begin{equation}\label{extenJ}\begin{split}
&\bar J_1\xi_2=\xi_3, \ \bar J_1\xi_3=-\xi_2, \\
&\bar J_1|_{\sfD}=J_1. 
\end{split}
\end{equation}
If $E\otimes \CC=
E^{1,0}\oplus E^{0,1}$ is the eigenspace decomposition for $\bar J_1$,
then $E^{1,0}=\langle \xi_2-i\xi_3\rangle\oplus {\sfD}^{1,0}$.
In order to prove Theorem \ref{hcr} below, we prepare the
following lemma, (the
proofs are essentially the same as those of 
\cite[Lemma 3.2,\,Proposition 3.3]{KA3}, \cite[Lemma 2.9,\, Section 2.1]{AK}.)

\begin{lemma}\label{hit} The following equality holds.
\begin{enumerate}
\item[(i)] For any $\mbi{x},\mbi{y}\in {\sfD}^{1,0}$, there is an element
 $\mbi{u}\in {\sfD}\otimes \CC$ such that 
$\displaystyle [\mbi{x},\mbi{y}]=a(\xi_2-i\xi_3)+\mbi{u}$
for some $a\in \RR$. Furthermore,
\item[(ii)] $\mbi{u}\in {\sfD}^{1,0}\otimes \CC$, that is
$J_1\mbi{u}=i\mbi{u}$.
\end{enumerate}
\end{lemma}

\begin{theorem}\label{hcr} 
Let $X$ be a simply connected non-compact
$qc$-manifold with the $qck$-distribution $V$.
Suppose $V$ generates $\RR^3\leq \Psh_{qc}(X)$. Then

$(1)$\, For each $\RR\leq \RR^3$, the quotient
$X/\RR$ is a complex contact manifold.
There is a holomorphic bundle over 
a hyperK\"ahler manifold $Y=X/\RR^3$:
\begin{equation}\label{1holbund}
\begin{CD}
 \CC@>>> X/\RR@>>> Y.
\end{CD}
\end{equation}  

$(2)$\ For each $\RR^2\leq \RR^3$,
the quotient manifold $X/\RR^2$ is
a strictly pseudoconvex $CR$-manifold which
has a pseudo-Hermitian $($Sasaki$)$
bundle:
\begin{equation}\label{PHSa}
\begin{CD}
\RR@>>> X/\RR^2@>>> Y.
\end{CD}\end{equation}
\vskip0.1cm
Let $M$ be a closed $qc$-manifold with the $qck$-group $T^3$.
Suppose $T^3$ lifts to an $\RR^3$-action to the universal covering 
of $M$. Then the following holds.

$(1)'$\ For each $S^1\leq T^3$,
the quotient oribifold $M/S^1$ 
supports a complex contact structure such that 
$M/S^1$ is the holomorphic orbibundle
over the hyperK\"ahler orbifold $Z=M/T^3$: 
\begin{equation}\label{holbund}
\begin{CD}
 T^1_\CC@>>> M/S^1@>>> Z.
\end{CD}
\end{equation}

$(2)'$\ For a torus $T^2\leq T^3$,
the quotient orbifold $M/T^2$
admits a strictly pseudoconvex $CR$-structure. Furthermore
this gives the pseudo-Hermitian\\ $($Sasaki$)$
orbibundle:
\begin{equation}\label{PHSa}
\begin{CD}
S^1@>>> M/T^2@>>> Z.
\end{CD}\end{equation}

The fundamental group $\Gamma=\pi_1(M)$
induces a nontrivial group extension:
$\displaystyle 1\ra\, \ZZ^2\ra\, \Gamma/\ZZ \lra\, Q\ra 1$.
In particular \eqref{holbund} is a nontrivial $T^1_\CC$-bundle.
$\Gamma$ also assigns to \eqref{PHSa} a \emph{nontrivial} group extension:
$\displaystyle 1\ra\, \ZZ\ra\, \Gamma/\ZZ^2 \lra\, Q\ra 1$.
\end{theorem}

\begin{proof}
Since $\RR^3$ acts properly on $X$, the orbit spaces 
$$X_1=X/\RR,\ X_2=X/\RR^2\ \mbox{and}\ X_3=X/\RR^3\,(=Y)$$
are smooth manifolds respectively for $\RR^k\leq \RR^3$ $(k=1,2,3)$.
For the $qc$-structure $(\eta, \{J_\al\}_{\al=1}^3)$ on $X$,
we prove $(1)$, $(1)'$.  
First, (i) of Lemma \ref{hit} implies 
$[{\sfD}^{1,0},{\sfD}^{1,0}]\subset E^{1,0}$.
For $\mbi{u}\in {\sfD}^{1,0}$,
as $\xi_2$ generates $t_2\in \RR^3$,  
$J_1[\xi_2,\mbi{u}]=i[\xi_2,\mbi{u}]$. Similarly $J_1[\xi_3,\mbi{u}]=i[\xi_3,\mbi{u}]$.
Since $\displaystyle [\xi_2-i\xi_3, \mbi{u}]=[\xi_2,\mbi{u}]-i[\xi_3, \mbi{u}]\in 
\sfD\otimes \CC$, it follows
$\displaystyle [\xi_2-i\xi_3, \mbi{u}]\in {\sfD}^{1,0}$.
Noting \eqref{extenJ}, it follows
\begin{equation}\label{Eint}
[E^{1,0},E^{1,0}]\subset E^{1,0},
\end{equation} that is $\bar J_1$ is integrable on $E$.
The projection $p_1:X\ra X_1$ maps isomorphically $E$ onto $TX_1$
at each point of $X_1$, so $p_1$ induces an almost complex structure
$\hat J_1$ on $X_1$ such that $\displaystyle p_{1*} \bar {J_1}|_{E}=\hat J_1 p_{1*}$.
For $\cA,\cB\in E^{1,0}$, noting $[\cA,\cB]\in E^{1,0}$,
it follows  $p_{1*}[\cA,\cB]=
[p_{1*}\cA,p_{1*}\cB]\in TX_1^{1,0}$, that is $\hat J_1$
is a complex structure on $X_1$. 
Moreover, as $\RR^3$ leaves each $\xi_\al$ invariant,
both $p_{1*}\xi_2$ and $p_{1*}\xi_3$ induce the vector fields
$\hat\xi_2, \hat\xi_3$ respectively on $X_1$.
Then it follows 
$TX_1^{1,0}=p_{1*}E^{1,0}=\langle \hat\xi_2-i\hat\xi_3\rangle\oplus p_{1*}D^{1,0}$
so that
\begin{equation}\label{Dint}
[p_{1*}D^{1,0},p_{1*}D^{1,0}]=p_{1*}[D^{1,0},D^{1,0}]\equiv \langle
\hat\xi_2-i\hat\xi_3\rangle\ \,{\rm mod}\ \,p_{1*}D^{1,0}.
\end{equation}
By the definition,
$p_{1*}D^{1,0}$ is an (invariant) complex contact bundle on the complex manifold 
$(X_1,\hat J_1)$.

For $(1)'$, let $M=X/\Gamma$ (we use $X$ as the universal covering of $M$
with the same notations as those of $X$).
 There is the commutative diagram of principal bundles:
\begin{equation}\label{sasaki1}
\begin{CD}
1@>>>\ZZ @>>> \Gamma@>p_1>> \Gamma_1@>>> 1 \\ 
@. @VVV  @VVV @VVV @.\\
@.\RR@>>> X@>p_1>> X_1 @.\\
@.@VVV @VVV  @VVV @.\\
@.S^1@>>> M @>\nu>> X_1/\Gamma_1 @.
\end{CD}
\end{equation}such that $\ZZ=\Gamma\cap \RR$  
and the bottom sequence is an orbibundle.
Taking the quotient of $X_1$ by $\Gamma_1$,
$M/S^1=X_1/\Gamma_1$ is a \emph{complex contact orbifold} induced 
from that of $X_1$. 
Putting $\ZZ^3=\RR^3\cap \Gamma$, consider
the equivariant principal bundle
$\displaystyle (\ZZ^3,\RR^3)\ra\, (\Gamma,X)\stackrel{\tilde q}{\lra}\,(Q,X_3)$
whose quotient gives an orbibundle: $\displaystyle T^3\ra M\lra Z=X_3/Q$.
Recall $\displaystyle t^*\eta=\eta,\, \ t_*J_\al=J_\al t_*|_\sfD$ \, 
$({}^\forall\, t\in \RR^3,\ \al=1,2,3)$ (\cf Proposition \ref{soR3}).
Noting $d\eta(V,\mbi{u})=0, \, d\eta(V,V)=0$ $(\mbi{u}\in \sfD)$ from \eqref{ijk},
there is a K\"ahler form $\Om_\al$ such that $\Om=\Om_1i+\Om_2j+\Om_3k$ on $X_3$ 
satisfying (\cf \cite{KA2}):
\begin{equation}\label{hKaehler}
d\eta=\tilde q^*\Om\ \, \mbox{on}\ \, X.
\end{equation}
As $\tilde q_* : \sfD\ra TX_3$ is an isomorphism at each point of $X_3$,
there is an almost complex structure $\sf J_\al$ on $X_3$ such
that $\displaystyle \tilde q_*J_\al={\sf J_\al} \tilde q_*$ on $\sfD$.
It follows from \eqref{hKaehler} that
 $\displaystyle \Om_\al({\sf J_\al}\mbi{x},{\sf J_\al}\mbi{y})
=\Om_\al(\mbi{x},\mbi{y})$ $(\mbi{x},\mbi{y}\in TX_3)$. Moreover, 
$\displaystyle \Om_1({\sf J_1}\mbi{x},\mbi{y})=\Om_2({\sf J_2}\mbi{x},\mbi{y})=
\Om_3({\sf J_3}\mbi{x},\mbi{y})$.
Each $\sf J_\al$
turns to a complex structure on $X_3$ such that
$g=\Om_\al\circ \sf J_\al$ $(\al=1,2,3)$
is a hyperK\"ahler metric on the complex manifold $Y=X_3$.
As $\langle \hat\xi_2-i\hat\xi_3\rangle$ generates $\CC$,
the above complex contact structure on $(X_1,\hat J_1)$ admits a holomorphic 
principal bundle
$\displaystyle \CC\,\ra (X_1,\hat J_1)\stackrel{q}\lra (X_3, \sfJ_1)$
over the K\"ahler manifold $(X_3, \sf J_1)$.
The following commutative diagram of principal bundles 
associates the commutative diagram of group extensions:
\begin{equation}\label{twsasaki}
\begin{CD}
\RR &=& \RR@.  \\
\downarrow@. \downarrow@. @. \\
\RR^3@>>> X@>\tilde q>> X_3\\
@VVV @Vp_1 VV  ||@.\\
\RR^{2}@>>> X_1 @> q>> X_3\,,\\
\end{CD}
\end{equation}  
\begin{equation}\label{2twsasaki}
\begin{CD}
@.\ZZ &=& \ZZ @.  \\
@.\downarrow@. \downarrow@. @. \\
1@>>>\ZZ^3@>>> \Gamma @>\tilde q>> Q@>>> 1\\
@.@VVV @V p_1 VV  ||@.\\
1@>>>\ZZ^{2}@>>> \Gamma_1 @>q>> Q@>>> 1.\\
\end{CD}
\end{equation}Taking the quotients gives a holomorphic orbibundle:
\begin{equation*}
\begin{CD}
\RR^2/\ZZ^2@>>> X_1/\Gamma_1@>\hat q>> X_3/Q\\
 ||@. ||@. ||@.\\
T^1_\CC @>>> M/S^1@>\hat q >> Z. \\
\end{CD}
\end{equation*}
This proves $(1)$, $(1)'$.
Next let 
$\displaystyle (\ZZ^2,\RR^2)\ra\, (\Gamma, X)\stackrel{\mu}\lra\, (\Gamma_2, X_2)$
be the equivariant principal bundle
where $\RR^2=\langle \xi_2,\,\xi_3\rangle$.
$\mu_*$ maps $\sfD$ isomorphically to $\mu_*(\sfD)$.
Put $\displaystyle \mu_*(\sfD)=F$ which is the codimension $1$-subbundle of
$TX_2$. If $J'_1$
is the induced almost complex structure on $F$ such that $\mu_*J_1=
J'_1\mu_*$, then $\displaystyle \mu_*:({\sfD}^{1,0},J_1)\ra (F^{1,0},J'_1)$
 is an isomorphism.
Since $\displaystyle [{\sfD}^{1,0},{\sfD}^{1,0}]\subset \langle \xi_2-i\xi_3\rangle 
\oplus {\sfD}^{1,0}$ by Lemma \ref{hit},
it implies $\displaystyle \mu_*[{\sfD}^{1,0},{\sfD}^{1,0}]=
[F^{1,0},F^{1,0}]\subset F^{1,0}$, that is
$J'_1$ is integrable on $F$. 
As $\eta_1(\{\xi_2,\xi_2\})=0$, $\eta_1$ induces
a $1$-form $\hat \eta_1$ on $X_2$
such that $\mu^*\hat \eta_1=\eta_1$. Put $\hat\xi_1=\mu_*(\xi_1)$ so that
$\hat\eta_1(\hat\xi_1)=1$. Since
$\displaystyle \mu^*(\hat \eta_1\we (d\hat\eta_1)^{2n})=
\eta_1\we (d\eta_1)^{2n}|_{\{\xi_1,\sfD\}}\neq 0$,
it follows  $\displaystyle \hat \eta_1\we (d\hat\eta_1)^{2n}\neq 0$.
$\displaystyle (\hat\eta_1, J_1')$ is a strictly pseudoconvex
pseudo-Hermitian structure on $X_2$ with Reeb field $\hat\xi_1$.
The equivariant principal bundle: 
$\displaystyle (\ZZ,\RR)\ra\, (\Gamma_2, X_2)\lra\, (Q,X_3)$
gives rise to a pseudo-Hermitian Sasaki orbibundle:
\begin{equation}\label{PSSa}\begin{CD}
S^1@>>>X_2/\Gamma_2@>>> Z=X_3/Q.
\end{CD}\end{equation}
This proves $(2)$, $(2)'$. 
Finally taking into account the commutative diagram,
\begin{equation}\label{2cocycle}
\begin{CD}
@.\ZZ &=& \ZZ @.  \\
@.@VVV @VVV @. \\
1@>>>\ZZ^2@>>> \Gamma_1 @>>> Q@>>> 1\\
@.@Vp_2 VV @V p_2 VV  ||@.\\
1@>>>\ZZ @>>> \Gamma_2 @>q_1>> Q@>>> 1\, ,\\
\end{CD}
\end{equation}
the group extension in the middle 
gives a $2$-cocycle $[f]\in H^2(Q,\ZZ^2)$.
As the projection $p_2:\ZZ^2\ra \ZZ$ induces a $2$-cocycle
$[p_2(f)]\in H^2(Q,\ZZ)$ which represents the
group extension $\displaystyle 1\ra\,\ZZ\ra\, \Gamma_2 \stackrel{q_1}\lra Q\,\ra 1$.
Since this group extension is obtained from the pseudo-Hermitian Sasaki orbibundle
of \eqref{PSSa}, it does not split, that is $[p_2(f)]\neq 0$.
Thus $[f]\neq 0$, the above holomorphic bundle \eqref{holbund} is
a nontrivial \emph{smooth} bundle. 
\end{proof}

\begin{eg}   
The quaternionic Heisenberg Lie group $\cM$ 
has the $qck$-group $\RR^3=\RR i+\RR j+\RR k={\rm Im}\,\HH$
{\rm (\cf Section \ref{sub:Ca})}.
 Taking $\RR\,(=\RR i)$,
the quotient Lie group $\cL=\cM/\RR$ is a  
complex Lie group which admits a complex contact structure
{\rm (\cf \cite{KA4})}. A holomorphic fibering gives an exact sequence:
$\displaystyle 1\ra\, \CC\ra\, \cL\lra\, \HH^n\ra 1$.
$\cL$ is called \emph{Iwasawa complex Lie group}. Let
$\displaystyle \Sim(\cL)=\cL\rtimes (\Sp(n)\cdot S^1\times \RR^+)$
be the holomorphic subgroup preserving the complex contact structure of $\cL$.
A $(4n+2)$-manifold $M$ locally modeled on $\cL$ with coordinate changes 
lying in $\Sim(\cL)$ is said to a \emph{complex contact similarity manifold}.
The following is shown similarly as in the proof of \cite{KA4}. 
\begin{pro}\label{ccstruc}
Let $M$ be a compact complex contact similarity manifold.
If the developing map $\displaystyle \dev:
\tilde M\lra \cL$ is injective,
then $M$ is holomorphically isometric to 
the infranilmanifold $\cL/\Gamma$ 
$(\Gamma\leq \cL\rtimes (\Sp(n)\cdot S^1))$
or an infra-Hopf manifold 
$\displaystyle \cL-\{\mbi{0}\}/\Gamma$ $(\Gamma\leq 
\Sp(n)\cdot S^1\times \RR^+)$.
\end{pro}
\end{eg}

\section{The standard $qc$-structure $\sfD_0$ for $(\cM,\om_0)$}\label{sub:Ca}
Recall the $4n+3$-dimensional (positive definite)
quaternionic Heisenberg Lie group $\mathcal M$
from \cite{AK}.
Put $\mbi{t}=(t_1,t_2,t_3), \mbi{s}=(s_1,s_2,s_3)\in\RR^3=\mbox{Im}\HH$, and
$z={}^t(z_1,\ldots,z_n), w={}^t(w_1,\ldots,w_n)\in \HH^n$ and so on.
Then $\cM$ is the product $\RR^3\times \HH^n$ with group law:
\[(\mbi{t},z)\cdot (\mbi{s},w)=(\mbi{t}+\mbi{s}-\mbox{Im}\langle z,w\rangle,\,z+w).\]
where $\langle z,w\rangle={}^t\bar zw$
is the Hermitian inner product.
As $\cM$ is nilpotent such that the center is the commutator subgroup
 $[\mathcal M,\mathcal M]=\RR^3$ consisting of elements $(t,0)$.
Let $\om_0$ be an ${\rm Im}\HH$-valued $1$-form on $\mathcal M$ defined by
\begin{equation}\label{Im-form}\begin{split}
\om_0=dt_1i+dt_2j+dt_3k+{\rm Im}\langle z,dz\rangle.
\end{split}
\end{equation} 
Put $\om_0=\om_1i+\om_2j+\om_3k$.
Denote by $\displaystyle {\sfD}_0$
the codimension $3$-subbundle $\displaystyle {\rm ker}\, \om_0$ on $\cM$.
As in \eqref{hypercom}, three endomorphisms $\{J_\al\}_{\al=1}^3$
are defined by  
$\displaystyle J_\ga=(d\om_\be|{\sf D_0})^{-1}\circ
(d\om_{\al}|{\sf D_0}):{\sf D_0} \ra {\sf D_0}$
on $\sfD_0$ \,$((\alpha, \beta,\gamma)\sim(1,2,3))$.
For the projection $\pi:\cM\ra \HH^n$,
$\pi_*:\sfD_0\ra\, T\HH^n$ is
the isomorphism at each point of $\HH^n$.
Let $\{i,j,k\}$ be the standard quaternionic structure on $\HH^n$.
We define a quaternionic structure 
$\{\tilde J_\al\}_{\al=1}^3$ on $\sfD_0$ by
$\pi_*\tilde J_1\mbi{u}=(\pi_*\mbi{u})\bar{i}$,\
$\pi_*\tilde J_2\mbi{u}=(\pi_*\mbi{u})\bar{j}$,\
$\pi_*\tilde J_3\mbi{u}=(\pi_*\mbi{u})\bar{k}$ respectively.
A calculation shows
that the quaternionic structure 
$\{\tilde J_\al\}_{\al=1}^3$ on $\sfD_0$
 coincides with $\{J_\al\}_{\al=1}^3$,
that is $J_\al=\tilde J_\al$ $(\al=1,2,3)$. (Compare \cite[Section 6.2]{AK1}.) 
By the definition, the pair $(\sfD_0,\{J_\al\}_{\al=1}^3)$ is 
the standard $qc$-structure on $\cM$ with
the $qck$-distribution 
$\displaystyle V_0=\langle\,  \frac d{dt_1},\,
 \frac d{dt_2},\, \frac d{dt_3}\,\rangle$ generating the center $\RR^3$.
By Theorem A,
${\Aut}_{qc}(\cM)=\cM\rtimes({\rm Sp}(n)\cdot \Sp(1)\times \RR^+)$ 
is the full group of $qc$-automorphisms of $\mathcal M$.
Recall \cite{KA2} that the action of an element
$\al=((s,z),\lambda\cdot A\cdot a)\in\mathop{\Aut}_{qc}(\mathcal M)$ on
$\cM$ is obtained as $({}^\forall\,(t,w)\in \cM)$ :
\begin{equation}\label{qc-action}\begin{split}
\al(t,w)&= 
(s,z)(\lambda^2\cdot a\cdot t\cdot\bar a,\, \lam\cdot A\cdot w\cdot \bar a)\\
&=(s+\lambda^2 a\cdot t\cdot\bar a-
{\rm Im}\langle z,\,\lambda\cdot A\cdot w\cdot \bar a\rangle,\, z+
\lambda\cdot A\cdot w\cdot \bar a).
\end{split}
\end{equation}
 For the standard $qc$-form $\om_0$ (\cf \eqref{Im-form}), it follows
 
\begin{equation}\label{actionqc}
\al^*\om_0=\lambda^2 a\cdot \om_0\cdot\bar a.
\end{equation} 
Then $\displaystyle \Psh_{qc}(\cM, \om_0,\{J_\al\}_{\al=1}^3)=
\cM\rtimes({\rm Sp}(n)\cdot \Sp(1))$ by Definition \ref{almosth}.

\section{Complete hyperK\"ahler structure on $\RR\times X$}\label{producthyp}
If $X$ is a $qCR$-manifold, then it is known that the product
$\RR^+\times X$ with the cone metric
admits an (incomplete) hyperK\"ahler structure.
This construction is not applicable to $qc$-manifolds with 
the $qck$-distribution $\RR^3$. (See Note \ref{actionSpM}.) 
However we shall construct a hyperK\"ahler metric on the product
$\RR\times X$
suitable for a $qc$-manifold $X$ with $qck$-group $\RR^3$.
Let $\displaystyle (X,\sfD,\om,\{J_\al\}_{\al=1}^3)$ be a simply connected non-compact
$qc$-manifold with the $qck$-group $\RR^3$ where
$\om=\om_1i+\om_2j+\om_3k$ is an ${\rm Im}\,\HH$-valued $1$-form.
By  \eqref{hKaehler} of Theorem \ref{hcr}, there is a principal
bundle $\displaystyle  \RR^3\ra X\stackrel{\pi}\lra Y$
over a hyperK\"ahler manifold $Y$.
$\RR^3$ induces
the distribution $\displaystyle V=\langle \frac d{dt_1},\, \frac d{dt_2},\, 
\frac d{dt_3}\rangle$ as before.
As $\displaystyle \om_\al(\frac d{dt_\be})=\delta_{\al\be}$
with \eqref{Im-form}, 
each dual form satisfies
\begin{equation}\label{equiomt}
dt_\al\equiv \om_\al \ \, ({\rm mod}\, \sfD). 
\end{equation} 

Choose the coordinate $t_0$ for $\RR$ in
the product $\RR\times X$. Then
the canonical distribution $\RR^4$
is endowed with the standard quaternionic structure such as

\begin{equation}\label{qr4}
\begin{split}
&\bar J_\al \frac d{dt_\al}=\frac d{dt_0}\, \ (\al=1,2,3),\\
&\bar J_\al \frac d{dt_0}=-\frac d{dt_\al}\,.
\end{split}
\end{equation}
Replace $\RR^4$ by
$\displaystyle\HH=\langle \frac d{dt_0},\,\frac d{dt_1},\, \frac d{dt_2},\, 
\frac d{dt_3}\rangle$.
As there is a decomposition $$T(\RR\times X)=\RR\times \RR^3\oplus
\sfD=\HH\oplus \sfD,$$ extend a quaternionic structure $\{\bar J_\al\}$ on
$\displaystyle \HH\oplus \sfD$ to be
\begin{equation}\label{extenJH}
\bar J_\al=\begin{cases}
\bar J_\al &\mbox{on}\ \HH\\
J_\al & \mbox{on}\ \sfD.
\end{cases}
\end{equation}

The product group $\RR\times \Psh_{qc}(X)$ acts on
$\RR\times X$ as usual:
\[ (t,h)(s,x)=(t+s,hx)\ \ ({}^\forall\,(s,x)\in \RR\times X).\] 
There is also a principal bundle
$\displaystyle \HH\ra\, \RR\times X\stackrel{\tilde\pi}\lra\, Y$.
Denote the subgroup   
\begin{equation}\begin{split}\label{subhygroup}
G=\{h\in \Psh_{qc}(X)&\mid\ h^*(dt_1i+dt_2j+dt_3k)=\\
&a_h(dt_1i+dt_2j+dt_3k)\overline{a_h}\,, \  ({}^\exists\, a_h\in \Sp(1))\,\}.
\end{split}\end{equation}
The normal subgroup $\displaystyle \HH\leq \RR\times \Psh_{qc}(X)$ 
is obviously contained in $G$ (see {\bf (i)} of Note \ref{actionSpM}).

\begin{pro}\label{hyperKX}
Let $X$ be a simply connected non-compact $qc$-manifold with 
the $qck$-group $\RR^3$. Suppose $Y$ is complete.
Then there is a complete hyperK\"ahler metric $g_0$ on the
quaternionic manifold
$(\RR\times X,\, \{\bar J_\al\}_{\al=1}^3)$ 
such that the quaternionic isometry group 
$\Isom_{hK}(\RR\times X)$ is $\RR\times G$.
\end{pro}
 
\begin{proof}
According to the decomposition, an arbitrary element 
is described as $\mbi{a}+\mbi{u}$ $(\mbi{a}\in \HH,
\mbi{u}\in \sfD)$. 
Define a Riemannian metric on $\RR\times X$ to be

\begin{equation*}\begin{split}\label{hyRiemann}
g_0(\mbi{a}+\mbi{u},\mbi{b}+\mbi{v})&=
dt_0(\mbi{a})\cdot dt_0(\mbi{b})+
dt_1(\mbi{a})\cdot dt_1(\mbi{b})+dt_2(\mbi{a})\cdot dt_2(\mbi{b})\\
&\ \ \ \ \ + dt_3(\mbi{a})\cdot dt_3(\mbi{b})+d\om_1(\bar J_1\mbi{u},\mbi{v}).
\end{split}\end{equation*}
Note 
$\displaystyle d\om_\al(\frac d{dt_0},\HH\oplus \sfD)=0$ by the decomposition.
Since $\displaystyle V=\langle \frac d{dt_1},\frac d{dt_2},\frac d{dt_3}\rangle$
is $\RR^3$, $\displaystyle d\om_\al(\frac d{dt_i},\frac d{dt_j})=0$ $(i,j=1,2,3)$
and $d\om_\al(V,\sfD)=0$ (\cf \eqref{ijk}). 
Using the projection
$\displaystyle p:T(\RR\times X)=\HH\oplus \sfD
\lra \HH$  on the first factor, it follows 

\begin{equation}\label{hypKRiemann}\begin{split}
g_0(\mbi{a}+\mbi{u},\mbi{b}+\mbi{v})=\bigl(p^*
(\mathop{\sum}_{i=0}^3{dt_i\cdot dt_i})+&
d\om_1\circ \bar J_1\bigr)(\mbi{a}+\mbi{u},\mbi{b}+\mbi{v})\\  
&\ \  (\mbi{a},\mbi{b}\in \HH,\ \mbi{u},\mbi{v}\in \sfD).
\end{split}\end{equation}
Since $\bar J_\al$ preserves the decomposition which also
 leaves invariant $\displaystyle \mathop{\sum}_{i=0}^3{dt_i\cdot dt_i}$,\,
 $d\om_\al\circ J_\al$ respectively, we obtain
\begin{equation*}\begin{split}
g_0(\bar J_\al(\mbi{a}+\mbi{u}),\bar J_\al(\mbi{b}+\mbi{v}))&=
g_0(\bar J_\al\mbi{a}+J_\al(\mbi{u}),\bar J_\al(\mbi{b})+J_\al(\mbi{v}))\\
&=g_0(\mbi{a}+\mbi{u},\mbi{b}+\mbi{v}).
\end{split}\end{equation*}
In order to prove $g_0$ is a hyperK\"ahler metric on $\RR\times X$,
put
\begin{equation}\label{2closesform}
\Theta_\al(\mbi{x},\mbi{y})=g_0(\mbi{x},\bar J_\al \mbi{y})
\ \ (\mbi{x},\mbi{y}\in T(\RR\times X),\, \al=1,2,3).
\end{equation}
It is easy to check that
\begin{equation*}\begin{split}
dt_0\circ \bar J_\al &=dt_\al\ (\al=1,2,3),\\
dt_1\circ \bar J_1&=-dt_0,\ dt_1\circ \bar J_2=dt_3,\ dt_1\circ \bar J_3=-dt_2,\\
dt_2\circ \bar J_1&=-dt_3,\ dt_2\circ \bar J_2=-dt_0,\ dt_2\circ \bar J_3=dt_1,\\
dt_3\circ \bar J_1&=dt_2,\ dt_3\circ \bar J_2=-dt_1,\ dt_3\circ \bar J_3=-dt_0.\\
\end{split}\end{equation*}
A calculation shows 
\begin{equation*}\begin{split}
\Theta_1(\mbi{a}+\mbi{u},\mbi{b}+\mbi{v})&=
g_0(\mbi{a}+\mbi{u},\bar J_1\mbi{b}+J_1\mbi{v})\\
&=\sum_{i=0}^3dt_i(\mbi{a})dt_i(\bar J_1\mbi{b})
+d\om_1(\bar J_1\mbi{u},\bar J_1\mbi{v})\\
&=2(dt_0\we dt_1+dt_3\we dt_2)(\mbi{a},\mbi{b})+d\om_1(\mbi{u},\mbi{v})\\
&=\bigl(2p^*(dt_0\we dt_1+dt_3\we dt_2)+d\om_1\bigr)(\mbi{a}+\mbi{u},\mbi{b}+\mbi{v}),\\
\Theta_2(\mbi{a}+\mbi{u},\mbi{b}+\mbi{v})&=
\sum_{i=0}^3dt_i(\mbi{a})dt_i(\bar J_2\mbi{b})
+d\om_2(\bar J_2\mbi{u}, \bar J_2\mbi{v})\\
&=\bigl(2p^*(dt_0\we dt_2+dt_1\we dt_3)
+d\om_2\bigr)(\mbi{a}+\mbi{u},\mbi{b}+\mbi{v})\\
\Theta_3(\mbi{a}+\mbi{u},\mbi{b}+\mbi{v})&=
\sum_{i=0}^3dt_i(\mbi{a})dt_i(\bar J_3\mbi{b})
+d\om_3(\bar J_3\mbi{u}, \bar J_3\mbi{v})\\
&=\bigl(2p^*(dt_0\we dt_3+dt_2\we dt_1)
+d\om_3\bigr)(\mbi{a}+\mbi{u},\mbi{b}+\mbi{v}).
\end{split}\end{equation*}

\begin{equation}\label{formof 2}
\Theta_\al=2p^*(dt_0\we dt_\al+dt_\ga\we dt_\be)+d\om_\al\ \
((\al,\be,\ga)\sim (1,2,3)).
\end{equation}
Each $\Theta_\al$ is closed $(\al=1,2,3)$ and
$\displaystyle \Theta_\al(\bar J_\al\mbi{x},\mbi{y})=
g_0(\mbi{x},\mbi{y})$ by \eqref{2closesform}. Thus
$g_0$ is a hyperK\"ahler metric on $\RR\times X$.
In particular, the principal bundle
\begin{equation}\label{metrichk}
 \HH\ra\, (\RR\times X,g_0)\stackrel{\tilde\pi}\lra\, (Y,\hat g_0)
\end{equation}
turns to a Riemannian submersion where $\hat g_0$ is the induced hyperK\"ahler
metric from $d\om_\al\circ J_\al$ (\cf \eqref{hKaehler}).
If $(Y,\hat g_0)$ is complete, then $(\RR\times X,g)$ will be also complete
by the construction.
Finally in view of \eqref{subhygroup},
the subgroup $\displaystyle \RR\times G$
is the quaternionic isometry group $\displaystyle\Isom_{hK}(\RR\times X,g_0)$. 
\end{proof}

\begin{theorem}\label{hyperM}
There is a complete hyperK\"ahler metric $g_0$ on the quaternionic space
$\HH^{n+1}$ $(n\geq 1)$ such that the quaternionic isometry group is
$\Isom_{hK}(\HH^{n+1},g_0)=\HH\rtimes (\Sp(n)\cdot \Sp(1))$. 
In particular, $g_0$ is not equivalent to the standard quaternionic metric
up to a quaternionic isometry. 
\end{theorem}

\begin{proof}
Take $\cM$ with $qck$-group $\RR^3=C(\cM)$.
We apply Proposition \ref{hyperKX} to $\RR\times \cM$.
In the decomposition $T(\RR\times \cM)=\HH\oplus \sfD_0$,
the quaternionic structure $\{J_\al\}_{\al=1}^3$ on $\sfD_0$
maps quaternionically onto the standard quaternionic structure on $\HH^n$ (\cf Section \ref{sub:Ca}). Noting \eqref{extenJ},
the quaternionic structure $(\RR\times \cM,\,\{\bar J_\al\}_{\al=1}^3)$
is equivalent with the standard quaternionic structure $(\HH^{n+1},\{i,j,k\})$ 
where an identification is 
$\displaystyle \RR\times \cM=\RR\times \RR^3\times \HH^n=\HH\times \HH^n=\HH^{n+1}$.
By the complete K\"ahler metric $\hat g_0$
on the base space $\HH^n$, using
the projection $\tilde \pi$ of \eqref{metrichk},
$(\HH^{n+1},g_0)$ will be a complete hyperK\"ahler manifold.
As $\displaystyle \Psh_{qc}(\cM,\om_0)=\cM\rtimes (\Sp(n)\cdot \Sp(1))$ 
by \eqref{actionqc},
the subgroup $G$ preserving both $\sum_{i=1}^3 dt_i\cdot dt_i$
and $d\om_1\circ J_1$
is isomorphic to $\displaystyle \RR^3\rtimes (\Sp(n)\cdot \Sp(1))$ 
$($see Note \ref{actionSpM}$)$.
Hence
$\displaystyle \Isom_{qK}(\HH^{n+1},g_0)=\HH\rtimes (\Sp(n)\cdot \Sp(1))$.
Since 
$\displaystyle \Isom_{qK}(\HH^{n+1},g_\HH)=\HH^{n+1}\rtimes (\Sp(n)\cdot \Sp(1))$
for the standard metric $g_\HH$, 
the hyperK\"ahler metric $g_0$ is not equivalent with
$g_\HH$ on $\HH^{n+1}$.
\end{proof}

\begin{note}\label{actionSpM}
{\rm {\bf (i)}}\,Put $\displaystyle d\mbi{t}=dt_1i+dt_2j+dt_3k$ for brevity.
Then $\displaystyle {\rm Re}\,\langle d\mbi{t}, d\mbi{t}\rangle=
dt_1\cdot dt_1+dt_2\cdot dt_2+dt_3\cdot dt_3$.
The conjugate of an element $a\in\Sp(1)$ leaves 
$\displaystyle {\rm Re}\,\langle d\mbi{t}, d\mbi{t}\rangle$ invariant, while
$\Sp(n)$ acts trivially on $d\mbi{t}$ (see \eqref{qc-action}). In fact,
$\displaystyle {\rm Re}\,\langle a\cdot d\mbi{t}\cdot\bar a
,\,a\cdot d\mbi{t}\cdot\bar a\rangle=
{\rm Re}\,(a\overline{d\mbi{t}}\bar a\cdot a d\mbi{t}\bar a)=
{\rm Re}\,(\overline{d\mbi{t}}\cdot{d\mbi{t}})=
{\rm Re}\,\langle d\mbi{t},d\mbi{t}\rangle$.
By \eqref{qc-action}, the elements of the form 
$(\mbi{0},z)\in((0,0,0),\HH^n)\leq \cM$
act on $((t_1,t_2,t_3),w)\in\cM$ as 
$\displaystyle(\mbi{0},z)((t_1,t_2,t_3),w)=((t_1,t_2,t_3)-{\rm Im}\langle
z,w\rangle, z+w)$, which 
do not preserve $dt_1\cdot dt_1+dt_2\cdot dt_2+dt_3\cdot dt_3$.
\noindent{\rm {\bf (ii)}}\, 
We shall
explain that the aforementioned cone-construction does not work.
Given a non-compact simply connected $qc$-manifold 
$\displaystyle (X,g,\om, \{J_\al\}_{\al=1}^3)$ 
with the $qck$-group $\RR^3$, $g(\mbi{x},\mbi{y})=
\sum_{i=1}^3\om_i(\mbi{x})\cdot \om_i(\mbi{y})+
d\om_1(J_1\mbi{x},\mbi{y})$ is the Riemannian metric 
invariant under $\Psh_{qc}(X,\om, \{J_\al\}_{\al=1}^3)$ on $X$ 
\rm {(\cf \eqref{Riemann})}.
Recall $\displaystyle g'=dt^2+t^2g$ is the cone metric on $\RR^+\times X$
$(t\in \RR^+)$ $(${\rm \cf} \cite{AK1}$)$.
The exact two-form  for $g'$ is
$\Om'_\al=d(t^2\cdot \om_\al)=2t dt\we \om_\al+t^2d\om_\al$.
Noting $\displaystyle V=\{ d/dt_1, d/dt_2, d/dt_3\}$ is
the $qck$-distribution generating $\RR^3$,
the quaternionic structure on $\RR^+\times X$ is defined by
$\displaystyle \bar J_\al d/dt_\al=t d/dt,\, \
\bar J_\al d/dt=- d/dt_\al$\, $(\al=1,2,3)$.
In particular
$\displaystyle\bar J_\al d/dt_\be= d/dt_\ga$.
If $g'$ happened to be a hyperK\"ahler metric, then 
it would satisfy
$\displaystyle \Om'_\al(\bar J_\al\mbi{x},\mbi{y})=g'(\mbi{x},\mbi{y})$
\ $(\al=1,2,3)$. 
However taking $\displaystyle \mbi{x}=\mbi{y}=d/dt_\be$,
it follows $\displaystyle g'(\mbi{x},\mbi{y})=(dt^2+t^2g)(d/dt_\be, d/dt_\be)=
t^2g(d/dt_\be,d/dt_\be)=t^2$, while
$\displaystyle \Om'_\al(\bar J_\al\mbi{x},\mbi{y})
=2t dt\we \om_\al(d/dt_\ga, d/dt_\be)
+t^2d\om_\al(d/dt_\ga, d/dt_\be)=0$.
The cone-metric is not hyperK\"ahler for
any $qc$-manifold $\displaystyle(X,g,\om, \{J_\al\}_{\al=1}^3,\RR^3)$.
\noindent{\rm {\bf (iii)}} Given a non-compact simply connected $qc$-manifold 
$\displaystyle (X,g,\om, \{J_\al\}_{\al=1}^3)$ 
with $V=\RR^3$ $($\cf Theorem $\ref{hyperM})$, there is
a Riemannian metric 
 $\displaystyle g_1=dt_0\cdot dt_0+\om_1\cdot \om_1+\om_2\cdot \om_2+
\om_3\cdot \om_3+d\om_1\circ J_1$ on $\RR\times X$.
Defining a quaternionic structure on $\RR\times X$
as in \eqref{qr4}, \eqref{extenJ}, we obtain $2$-forms
$\displaystyle \Om_\al=-4(\om_\al\we dt_0+\om_\be\we \om_\ga)+d\om_\al$
$((\al,\be,\ga)\sim (1,2,3))$ by a calculation, so $\Om_\al$ is not closed. However,
$\Isom_{hk}(\RR\times X,g_1)=\RR\times\Psh_{qc}(X)$ is the isometry group
which may act transitively on $\RR\times X$.
\end{note} 

\section{$\Psh_{qc}(X)$ with the abelian $qck$-distribution} \label{abelnormal} 
Let $\Psh_{qc}(X)=\Psh_{qc}(X,(\eta,\{J_\al\}_{\al=1}^3))$ be
the $qc$-Hermitian group defined in Definition \ref{almosth}.
When the $qck$-distribution $V$ generates $\RR^3$,
we observe that $\Psh_{qc}(X)$ is determined
exactly by using the Boothby-Wang fibering. (Compare \cite[Proposition 3.4]{OK}
for the case of pseudo-Hermitian Sasaki manifolds.)
As $\displaystyle N_{\Psh_{qc}(X)}(\RR^3)=\Psh_{qc}(X)$
by Corollary \ref{normalizer},
there is a principal bundle :
\begin{equation}\label{BW3}
\begin{CD}
\RR^3@>>> X@>p>> Y\\
\end{CD}
\end{equation}for which $p^*\Om=d\eta$ (\cf \eqref{hKaehler}). 
As in Definition \ref{almosth}, we introduce
\emph{hyperK\"ahler isometry group} as

\begin{definition}\label{hygr}
\begin{equation*}\label{qcIso}
\Isom_{hk}(Y, \Om,\{{\sf J_\al}\}_{\al=1}^3)=\bigl\{
h\in\Diff(Y)
\, \big|\  h^*\Om
=b \cdot\Om\cdot \bar b\,,\\
h_* \circ {\sf J}_k=\sum_{j=1}^{3}
b_{kj}{\sf J}_j\circ h_*\}
\end{equation*}for
$b\in C^\infty(Y, \Sp(1))$, and
$(b_{ij})\in C^\infty(Y,{\rm SO}(3))$ is its conjugate by $b$.
\end{definition}

Take $h\in \Psh_{qc}(X)$ such that
\begin{equation}\label{normala}
 h^*\eta=a\cdot \eta\cdot \bar a
\end{equation}
for some map $a\in C^\infty(X,\Sp(1))$ (\cf Definition \ref{almosth}).
For $t\in \RR^3$, $hth^{-1}\in \RR^3$ as above,
there induces a map $\hat h : Y\ra Y$ with the commutative diagram:
\begin{equation}\label{hatf}
\begin{CD}
X@>h>> X\\
@VpVV @VpVV\\
Y@>\hat h>> Y.\\
\end{CD}\end{equation}
Applying $t\in \RR^3$ to \eqref{normala},
$\displaystyle t^*a\cdot \eta \cdot \overline{t^*a}=t^*h^*\eta=h^*(hth^{-1})^*\eta=
h^*\eta=a\cdot \eta\cdot \bar a$,
it follows
$\displaystyle t^*a=\pm a\in C^\infty(X,\Sp(1))$ 
since $\eta(\,{}^t[\xi_1, \xi_2,\xi_3]\,)=
{}^t[\,i,j,k\,]$. Thus $t^*a=a$ and so $a$ induces a map 
$b\in C^\infty(Y,\Sp(1))$
such that $p^*b=a$.
Differentiate \eqref{normala}, then
$\displaystyle h^*d\eta=a\cdot d\eta\cdot \bar a$ on $\sfD$.
Since $\displaystyle h^*p^*\Om=a\cdot p^*\Om\cdot \bar a$ on $\sfD$,
the commutativity shows
$\displaystyle p^*{\hat h}^*\Om=p^*(b\cdot \Om\cdot \bar b)$
on $\sfD$, that is
\begin{equation}\label{downOM}
{\hat h}^*\Om=b\cdot \Om\cdot \bar b.
\end{equation}
As $h_*J_\al=\sum_{\be}a_{\al\be}J_\be h_*$ on $\sfD$,
it follows 
${\hat h}_*{\sf J}_\al=\sum_{\be}\, b_{\al\be}{\sf J}_\be {\hat h}_*$
on $Y$. Here $b_{\al\be}$ is induced from $b$.
By the definition, we have
${\hat h}\in \Isom_{hk}(Y, \Om,\{{\sf J}_\al\}_{\al=1}^3)$.
By \eqref{downOM}, let ${\hat h}^* \Om_\al=
\sum_{\be}\ b_{\al\be}\Om_\be$ where $b_{\al\be}\in C^\infty(Y,\SO(3))$.

\begin{lemma}\label{aconst}
Assume $\dim\,Y\geq 8$. 
Then $b_{\al\be}$ is constant, that is an element of $\Sp(1)$. 
In particular $\displaystyle h^*\eta=a\cdot \eta\cdot \bar a$
for an element $a\in\Sp(1)$. 
\end{lemma}

\begin{proof} Differentiate \eqref{downOM} such that
$0=d{\hat h}^* \Om_\al=
\sum_{\be}\ d b_{\al\be}\we \Om_\be$.
So we may put $\theta_1\we \Om_1+\theta_2\we \Om_2+\theta_3\we \Om_3=0$
where $\theta_k=d b_{\al k}$ for each $\al$ $(k=1,2,3)$.
It suffices to show $\theta_k=0$ on $Y$. Let $g=\Om_\al\circ {\sf J_\al}$.
Choose $\mbi{x}\in TY$ such that $g(\mbi{x},\,\mbi{x})=1$.
There is a decomposition :
$\displaystyle TY=\{\mbi{x},{\sf J_1}\mbi{x},{\sf J_2}\mbi{x},{\sf J}_3\mbi{x}\}\oplus
\{\mbi{x},{\sf J}_1\mbi{x},{\sf J}_2\mbi{x},{\sf J}_3\mbi{x}\}^{\perp}$.
\ {\bf Case 1}. Let $\mbi{y}\in \{\mbi{x},{\sf J}_1\mbi{x},{\sf J}_2\mbi{x},{\sf J}_3\mbi{x}\}^{\perp}$.
\begin{equation*}\begin{split}
\theta_1\we \Om_1(\mbi{y},{\sf J}_1\mbi{x},\mbi{x})=&
\theta_1(\mbi{y})\Om_1({\sf J}_1\mbi{x},\mbi{x})+
\theta_1(\mbi{x})\Om_1(\mbi{y},{\sf J}_1\mbi{x})+
\theta_1({\sf J}_1\mbi{x})\Om_1(\mbi{x},\mbi{y})\\
=&\theta_1(\mbi{y})g(\mbi{x},\mbi{x})-
\theta_1(\mbi{x})g(\mbi{x},\mbi{y})-
\theta_1({\sf J}_1\mbi{x})g({\sf J}_1\mbi{x},\mbi{y})
=\theta_1(\mbi{y}).\\
\theta_2\we \Om_2(\mbi{y},{\sf J}_1\mbi{x},\mbi{x})=&
\theta_2(\mbi{y})\Om_2({\sf J}_1\mbi{x},\mbi{x})+
\theta_2(\mbi{x})\Om_2(\mbi{y},{\sf J}_1\mbi{x})+
\theta_2({\sf J}_1\mbi{x})\Om_2(\mbi{x},\mbi{y})\\
=&-\theta_2(\mbi{y})\Om_3(\mbi{x},\mbi{x})-
\theta_2(\mbi{x})g({\sf J}_3\mbi{x},\mbi{y})-
\theta_2({\sf J}_1\mbi{x})g({\sf J}_2\mbi{x},\mbi{y})\\
=&0.\\
\theta_3\we \Om_3(\mbi{y},{\sf J}_1\mbi{x},\mbi{x})=&
\theta_3(\mbi{y})\Om_3({\sf J}_1\mbi{x},\mbi{x})+
\theta_3(\mbi{x})\Om_3(\mbi{y},{\sf J}_1\mbi{x})+
\theta_3({\sf J}_1\mbi{x})\Om_3(\mbi{x},\mbi{y})\\
=&\theta_3(\mbi{y})\Om_2(\mbi{x},\mbi{x})+
\theta_3(\mbi{x})g({\sf J}_2\mbi{x},\mbi{y})-
\theta_3({\sf J}_1\mbi{x})g({\sf J}_3\mbi{x},\mbi{y})
=0.
\end{split}
\end{equation*}The above equality shows $\theta_1(\mbi{y})=0$
for all $\mbi{y}\in \{\mbi{x},{\sf J}_1\mbi{x},{\sf J}_2\mbi{x},{\sf J}_3\mbi{x}\}^{\perp}$.\\
{\bf Case 2}. For $\mbi{x}\in \{\mbi{x},{\sf J}_1\mbi{x},{\sf J}_2\mbi{x},{\sf J}_3\mbi{x}\}$,
choose $\mbi{y}\in \{\mbi{x},{\sf J}_1\mbi{x},{\sf J}_2\mbi{x},{\sf J}_3\mbi{x}\}^{\perp}$
such that $g(\mbi{y},\mbi{y})=1$. Note
${\sf J}_1\mbi{y},{\sf J}_2\mbi{y},{\sf J}_3\mbi{y}\in \{\mbi{x},{\sf J}_1\mbi{x},{\sf J}_2\mbi{x},
{\sf J}_3\mbi{x}\}^{\perp}$.
Then
\begin{equation*}\begin{split}
\theta_1\we \Om_1(\mbi{x},{\sf J}_1\mbi{y},\mbi{y})&=
\theta_1(\mbi{x})\Om_1({\sf J}_1\mbi{y},\mbi{y})+
\theta_1(\mbi{y})\Om_1(\mbi{x},{\sf J}_1\mbi{y})+
\theta_1({\sf J}_1\mbi{y})\Om_1(\mbi{y},\mbi{x})\\
&=\theta_1(\mbi{x})g(\mbi{y},\mbi{y})-
\theta_1(\mbi{y})g(\mbi{y},\mbi{x})+
\theta_1({\sf J}_1\mbi{y})g({\sf J}_1\mbi{x},\mbi{y})\\
&=\theta_1(\mbi{x}).\\
\theta_2\we \Om_2(\mbi{x},{\sf J}_1\mbi{y},\mbi{y})&=
\theta_2(\mbi{x})\Om_2({\sf J}_1\mbi{y},\mbi{y})+
\theta_2(\mbi{y})\Om_2(\mbi{x},{\sf J}_1\mbi{y})+
\theta_2({\sf J}_1\mbi{y})\Om_2(\mbi{y},\mbi{x})\\
&=0.\\
\theta_3\we \Om_3(\mbi{x},{\sf J}_1\mbi{y},\mbi{y})&=
\theta_3(\mbi{x})\Om_3({\sf J}_1\mbi{y},\mbi{y})+
\theta_3(\mbi{y})\Om_3(\mbi{x},{\sf J}_1\mbi{y})+
\theta_3({\sf J}_1\mbi{y})\Om_3(\mbi{y},\mbi{x})\\
&=0.\\
\end{split}
\end{equation*} Thus $\theta_1(\mbi{x})=0$.
This calculation implies similarly 
$\displaystyle \theta_1({\sf J_1}\mbi{x})=\theta_1({\sf J}_2\mbi{x})=\theta_1({\sf J}_3\mbi{x})=0$.
Hence it follows $\theta_1=d b_{\al 1}=0$ on $Y$.
Thus $b_{\al 1}\in C^\infty(Y,\SO(3))$ is a constant.
The above argument shows also $\theta_2=\theta_3=0$, that is
$d b_{\al 2}=d b_{\al 3}=0$ on $Y$.
The matrix $(b_{\al\be})\in C^\infty(Y,\SO(3))$ is constant and so
is $b\in C^\infty(Y,\Sp(1))$. 
Since $a=p^*b$, $a\in C^\infty(X,\Sp(1))$ is also constant.
\end{proof}

\begin{corollary}\label{abeliabcase}
When the $qck$-distribution $V$ is $\RR^3$, 
the $qc$-Hermitian group $($respectively
hyperK\"ahler isometry group$)$ is described as follows:
\begin{equation*}
\begin{split}
\Psh_{qc}(X, \om,\{J_\al\}_{\al=1}^3)=\bigl\{
h\in\Diff(X)&
\, \big|\  h^*\om
=a \cdot\om\cdot \bar a\, \ (a\in\Sp(1)), \\
 \ \ \ \ h_* \circ J_k=\sum_{j=1}^{3}
&a_{kj}J_j\circ h_*|_{\sfD},\ (a_{ij})\in{\rm SO}(3) \},\\
& \\
\Isom_{hK}(Y, \Om,\{{\sf J}_\al\}_{\al=1}^3)=\bigl\{\hat h\in \Diff(Y)&
\, \big|\  \hat h^*\Om =b \cdot\Om\cdot \bar b\,\ (b\in \Sp(1)), \\
\ \ \ \ \ \ \ \hat h_* \circ \sf J_k=\sum_{j=1}^{3}&
b_{kj}{\sf J}_j\circ \hat h_*\,, \ (b_{ij})\in{\rm SO}(3) \}.\
\end{split}\end{equation*}
\end{corollary}

\begin{pro}\label{exact3}
Suppose $H^1(Y;\RR^3)=0$.
There associates a natural exact sequence 
 with the principal bundle
$\displaystyle \RR^3\ra X\stackrel{p}\lra Y$ of \eqref{BW3}:
\begin{equation}\begin{CD}\label{eq:const}
1@>>> \RR^3@>>>\Psh_{qc}(X)@>\phi>>\Isom_{hK}(Y)@>>>1\,
\end{CD}\end{equation}
\end{pro}

\begin{proof}
The proof is almost similar to the argument of Section $4$ of \cite{OK}.
As usual put 
$\displaystyle \eta=\eta_1i+\eta_2j+\eta_3k,\ \Om=\Om_1i+\Om_2j+\Om_3k$. 
Then we have the equality $p^*\Om=d\eta$.
Suppose $V=\RR^3$ consists of one-parameter groups $\{\varphi^{(1)}_{t},
\varphi^{(2)}_{t},\varphi^{(3)}_{t};\, t\in \RR\}$.
If we choose a section $s:Y\ra X$ $(p \circ s = \mathrm{id}_{Y})$ such that
$X$ is equivalent with the trivial bundle $\RR^3\times Y$ by
a bundle map $f:\RR^3\times Y\ra X$ such as
$\displaystyle f((t_1,t_2,t_3),y)=\varphi^{(1)}_{t_1}\circ
\varphi^{(2)}_{t_2}\circ \varphi^{(3)}_{t_3}( s(y))$.
This gives a commutative diagram:
\begin{equation}\label{eq:commuta}
\begin{CD}
 @. \RR^3\times Y@>f>> \;  \;  \; X  \\
@.     {\rm pr}\searrow @.   
\! \! \swarrow {p} @ .   \\
@.   &  Y \  @. \\
\end{CD} 
\end{equation} Define a $1$-form $\theta=s^*\eta$ which satisfies
\begin{equation}\label{theta1}
d\theta=\Om.
\end{equation}
Let $\eta_0=\sum_{i=1}^3 dt_i +{\rm pr}^*\theta$
be the $1$-form on the product $\RR^3\times Y$ satisfying $d\eta_0={\rm pr}^*\Om$
compatible with the regular hyperK\"ahler structure. (See \cite[Proposition 3.1]{OK}.)
$\sfD_0=\ker\, \eta_0$ admits a quaternionic structure $\{J_\al\}_{\al=1}^3$
which is the pullback of the quaternionic structure $\{\sf J_\al\}_{\al=1}^3$
on $Y$. (See (3.5), (3.6) of \cite{OK}.)
Then it is easy to check that 
\begin{equation}\label{fform}
 f^*\eta=\eta_0,\
 f_*J_\al=J_\al f_*|_{\sfD_0}.
\end{equation} 
Then the $qc$-structure 
$\displaystyle (\RR^3\times Y,\,\eta_0,\,\{J_\al\}_{\al=1}^3,\, \{d/{dt_i}\}_{i=1}^3)$ is equivalent through the bundle map $f$ with 
$\displaystyle (X,\,\eta,\,\{J_\al\}_{\al=1}^3,\{\xi_i\}_{i=1}^3)$.
(Compare \cite[Proposition 3.1]{OK}.)
By the commutative diagram \eqref{hatf},
define $\phi(h)=\hat h$ for $h \in \Psh_{qc}(X)$.
Then $\phi: \Psh_{qc}(X)\ra \Isom_{hK}(Y)$ is a homomorphism.
In order to prove 
the exactness of \eqref{eq:const}, take $\hat h\in \Isom_{hk}(Y)$.
By Corollary \ref{abeliabcase}, $\hat h^*\Om=b\cdot \Om \cdot \bar b$\, 
for some $b\in \Sp(1)$.
Define $h_1:\RR^3\times Y\ra\, \RR^3\times Y$ to be
$\displaystyle h_1(\mbi{t},y)=(\mbi{t},\hat h(y))$ where $\mbi{t}=(t_1,t_2,t_3)\in\RR^3$.
Put
\begin{equation}\label{hgauge}
\eta'=\bar b\cdot h_1^*\eta_0 \cdot b.
\end{equation} By a calculation,
$\displaystyle dh_1^*\eta_0 =h_1^*{\rm pr}^*\Om
={\rm pr}^*\hat h^*\Om=b\cdot {\rm pr}^*\Om \cdot \bar b$, it follows
\begin{equation}\label{compatom}
d\eta'={\rm pr}^*\Om.
\end{equation}
Letting the obvious section $s':Y\ra\, \RR^3\times Y$, put $\theta'=(s')^*\eta'$
on $Y$. As $s'\circ {\rm pr}|_{0\times Y}={\rm id}$,
it follows ${\rm pr}^*\theta'|_{0\times Y}=\eta'|_{0\times Y}$.
Noting 
$$\eta'(\frac d{dt_1})=\bar b\cdot i\cdot b,\
\eta'(\frac d{dt_2})=\bar b\cdot j\cdot b,\
\eta'(\frac d{dt_3})=\bar b\cdot k\cdot b,$$
it implies 
$\displaystyle \eta'=\bar b\cdot d\mbi{t}\cdot b+{\rm pr}^*\theta'$
on $\RR^3\times Y$.
(Here $\displaystyle d\mbi{t}=dt_1 i+dt_2 j+dt_3 k$.)
Note by \eqref{compatom} 
$\displaystyle d\theta'=(s')^*d\eta'=
(s')^*{\rm pr}^*\Omega=\Omega$.
Since $\displaystyle d\theta'=\Om=d\theta$ on $Y$ from \eqref{theta1},
 $[\theta-\theta']\in H^1(Y,\RR^3)$.
By the hypothesis $H^1(Y;\RR^3)=0$,
there is a map $\lam:Y\ra \RR^3$ such that 
\begin{equation}\label{cocy1}
\theta-\theta'=d\lam.
\end{equation} 
Let $G:\RR^3\times Y\ra\, \RR^3\times Y$
be a gauge transformation defined by
 $\displaystyle G(\mbi{t},\,y)=(b\cdot(\mbi{t}+\lam(y))\cdot \bar b,\,y)$.
Noting ${\rm pr}\circ G={\rm pr}$,
\begin{equation*}\begin{split}
G^*\eta'&=G^*(\bar b\cdot d\mbi{t}\cdot b+{\rm pr}^*\theta')
=\bar b\cdot dG^*\mbi{t}\cdot b+{\rm pr}^*\theta' \\
&=\bar b\cdot d(b(\mbi{t}+\lam(y))\bar b)\cdot b+{\rm pr}^*\theta'
=d(\mbi{t}+{\rm pr}^*\lam)+{\rm pr}^*\theta'\\
&=d\mbi{t}+{\rm pr}^*\theta =\eta_0\ \, (\eqref{cocy1}).\\
\end{split}\end{equation*}
Put $h'=h_1\circ G :\RR^3\times Y\ra \RR^3\times Y$.
In fact, $\displaystyle h'(\mbi{t},y)
=(b\cdot (\mbi{t}+\lam(y))\cdot\bar b,\hat h(y))$.
Noting $\displaystyle h_1^*\eta_0=b\cdot \eta'\cdot \bar b$ by \eqref{hgauge},
a calculation shows
$\displaystyle h^*\eta_0
=b\cdot \eta_0\cdot \bar b$.
Since $f^*\eta=\eta_0$, 
we obtain $\displaystyle (fh' f^{-1})^*\eta=
b\cdot \eta\cdot \bar b$.
Put $\displaystyle {h}=f h' f^{-1}$ and so 
$\displaystyle p{h}\circ f=pf\circ h'$.
By \eqref{eq:commuta}, it follows
$\displaystyle p\circ {h}=\hat h\circ p$.
As
$\hat h_*{\sf J}_\al=\sum_{\be=1}^3 b_{\al\be}{\sf J}_\be \hat h_*$ 
and $\displaystyle p_*J_\al={\sf J}_\al p_*$,
we have
$\displaystyle p_*{h}_* J_\al= 
\sum b_{\al\be}{\sf J}_\be \hat h_* p_*
=p_*(\sum b_{\al\be} J_\be {h}_*):\sfD\ra\, TY$,
thus 
${h}_* J_\al=\sum_{\be=1}^3 b_{\al\be} J_\be {h}_*$.
This implies ${h}\in \Psh_{qc}(X)$ such that $\phi(h)=\hat h$ by \eqref{hatf}.
This shows $\phi$ is surjective.
Let $\phi(h)=\hat h$ as above.
Put $\displaystyle h_0=f^{-1}\circ h\circ f : \RR^3\times Y\ra \RR^3\times Y$. 
As in the argument of surjectivity, there is a map $\lam:Y\ra \RR^3$ 
such that $\displaystyle
h_0(\mbi{t},y)=(h_0\mbi{t}h_0^{-1}+\lam({\hat h}(y)),\hat h(y))$.
If $h^*\eta=a \cdot\eta\cdot \bar a$ $(a\in \Sp(1))$, then note
$\hat h^*\Om=a \cdot\Om\cdot \bar a$ (\cf \eqref{downOM}). 
Suppose $\phi(h)=1$. Then $\Om=a \cdot\Om\cdot \bar a$.
Noting each $\Om_i$ is a K\"ahler form for
$\Om=\Om_1i+\Om_2j+\Om_3k$, this equation implies $a=\pm 1$
so that $\displaystyle h^*\eta=\eta$. 
Since $f^*\eta=\eta_0$ by \eqref{fform},
this shows $\displaystyle h_0^*\eta_0=\eta_0$.
Moreover $h_0(\mbi{0},y)=(\lam(y), y)$ as above.
For any $(0,\mbi{w}_y)\in (0,TY)\subset T\RR^3\times TY$, it follows
$(h_0)_*(\mbi{w}_y)=(\lam_*(\mbi{w}_y)+\mbi{w}_y)$. A calculation shows
$\displaystyle \eta_0(\mbi{w}_y)=h_0^*\eta_0(\mbi{w}_y)=
\eta_0(\lam_*(\mbi{w}_y)+\mbi{w}_y)=
\eta_0(\lam_*(\mbi{w}_y))+\eta_0(\mbi{w}_y)$.
Thus $\eta_0(\lam_*(\mbi{w}_y))=0$, so $\lam_*(\mbi{w}_y)\in \sfD_0=\ker\,\eta_0$.
Since $\lam_*(\mbi{w}_y)\in T\RR^3=V$, it follows $\lam_*(\mbi{w}_y)=0$,
that is $\lam_*(TY)=0$. Hence $\lam$ is a constant $\mbi{s}\in \RR^3$.
Then $\displaystyle h_0(\mbi{t},y)=(h_0\mbi{t}h_0^{-1}+\mbi{s}, y)$.
As $\displaystyle \eta_0=d\mbi{t}+{\rm pr}^*\theta$,
\begin{equation*}\begin{split}
h_0^*\eta_0&=dh_0^*\mbi{t}+h_0^*{\rm pr}^*\theta
=d(h_0\mbi{t}h_0^{-1}+\mbi{s})+{\rm pr}^*\theta
       =d(h_0\mbi{t}h_0^{-1})+{\rm pr}^*\theta.\\
\end{split}\end{equation*}
Noting $h_0^*\eta_0=\eta_0$, $\displaystyle h_0\mbi{t}h_0^{-1}-
\mbi{t}=\mbi{c}$ for some constant $\mbi{c}\in \RR^3$.
Taking $\mbi{t}=\mbi{0}$, $\mbi{c}=\mbi{0}$ and so
$\displaystyle h_0\mbi{t}h_0^{-1}=\mbi{t}$.
Hence $h_0(\mbi{t},y)=(\mbi{t}+\mbi{s},y)=\mbi{s}(\mbi{t},y)$ so that
 $h_0=\mbi{s}\in \RR^3$. Since $f$ is $\RR^3$-equivariant,
the equation $\displaystyle h=f\circ h_0\circ f^{-1}$ implies $h=\mbi{s}$.
This proves the exactness of $\phi$.
\end{proof}

\section{Spherical homogeneous $qc$-manifolds}\label{Example1}

We determine spherical (uniformizable) homogeneous $qc$-manifolds 
following the method of spherical homogeneous $CR$-manifolds \cite{BS}.
 Let $\displaystyle (\rho, \dev): (\Aut_{qc}(M), M) 
\ra (\Aut_{qc}(S^{4n+3}), S^{4n+3})$
be a developing pair. Suppose that $M$ is homogeneous by a
subgroup $\cG$ of $\Aut_{qc}(M)$.
Put $G=\rho(\cG)\leq \Aut_{qc}(S^{4n+3})=\PSp(n+1,1)$.
We may assume $G$ is a \emph{non-compact} closed subgroup
taking the closure if necessary. 
Then $G\cdot p\,(=\dev(M))$ is a homogeneous domain of $S^{4n+3}$ for some $p\in \dev(M)$.
For a closed submanifold $L\subset S^{4n+3}$,
denote by $\mathop{\Aut}_{qc}(S^{4n+3}-L)$
the subgroup of ${\PSp}(n+1,1)$ whose elements leave $L$ invariant.
If $G$ has the radical,
then $G$ belongs to the maximal amenable Lie subgroup 
in $\PSp(n+1,1)$, that is $\Aut_{qc}(\cM)$ up to conjugate.
 Otherwise, $G$ is semisimple. 
It follows (\cite{CG}, \cite[Lemma 3.1]{KU})
 
\begin{lemma}\label{tran1}
If $G$ is non-compact semisimple but not ${\PSp}(n+1,1)$,
 then $G$ is one of the following groups up to conjugate.
In each case $G$ acts properly on $S^{4n+3}-L$ where $L$ is a sub-sphere.
\begin{enumerate}
\item[{\bf (1)}] $\displaystyle G=\Aut_{qc}(S^{4n+3}-S^{m-1})=
P({\rm O}(m,1)\cdot{\rm Sp}(1)
\times {\rm Sp}(n-m+1))$ $(1\leq m\leq n+1)$.
${\rm Sp}(n-m+1)$ is the maximal compact subgroup
which fixes $S^{m-1}\ (=\partial H^{m}_{\RR})$.

\item[{\bf (2)}] $\displaystyle G=\Aut_{qc}(S^{4n+3}-S^{2m-1})=
P({\rm U}(m,1)\cdot{\rm U}(1)\times {\rm Sp}(n-m+1))$  $(1\leq m\leq n+1)$.
${\rm Sp}(n-m+1))$ is the maximal compact subgroup which fixes 
$S^{2m-1}\ (=\partial H^{m}_{\CC})$.

\item[{\bf (3)}] $\displaystyle G=\Aut_{qc}(S^{4n+3}-S^{4m-1})=
{\rm Sp}(m,1)\cdot{\rm Sp}(n-m+1)$ $(1\leq m\leq n)$.
${\rm Sp}(n-m+1)$ is the maximal compact subgroup
which fixes $S^{4m-1}\ (=\partial H^{m}_{\HH})$.

\item[{\bf (4)}] $G={\Aut}_{qc}(S^{4n+3}-S^{2})=
{\rm SL}_2(\CC)\cdot{\rm Sp}(n)$.
${\rm Sp}(n)$ is the maximal compact subgroup
which fixes $S^2=\partial\HH_{{\rm Im}\,\HH}$.
In this case, it fixes $S^3=\partial \HH^1_{\HH}\supset 
S^2$ so this case reduces to case ${\bf (3)}$.
\end{enumerate}
\end{lemma}

\begin{pro}\label{tran2}
Let $G$ be a non-compact semisimple subgroup of $\PSp(n+1,1)$.
$(1)$ Suppose $G$ is not the whole group $\PSp(n+1,1)$.
Only $\displaystyle G=
{\rm Sp}(m,1)\cdot{\rm Sp}(n-m+1)$ acts transitively
on $S^{4n+3}-S^{4m-1}$ $(1\leq m\leq n)$.
The positive definite homogeneous $qc$-manifold
$S^{4n+3}-S^{4m-1}$ supports a principal bundle
$\displaystyle S^3\ra\, S^{4n+3}-S^{4m-1}\lra\, 
\HH^m_{\HH}\times \HH\PP^{n-m}$ in which $S^3$ is not a $qck$-group.
\  $(2)$ Suppose $G=\PSp(n+1,1)$. Then $\displaystyle S^{4n+3}=\PSp(n+1,1)/\Aut_{qc}(\cM)
=\Sp(n+1)\cdot \Sp(1)/\Sp(n)\cdot \Sp(1)$
in which $\Sp(1)$ is the $qck$-group of $\displaystyle \Sp(n+1)\cdot \Sp(1)=
\Psh_{qc}(S^{4n+3})=N_{\Psh_{qc}(S^{4n+3})}(\Sp(1))$.
\end{pro}

Suppose $G$ has the non-compact radical in $\PSp(n+1,1)$.
Then $G$ belongs to the maximal amenable Lie subgroup
$\displaystyle \cM\rtimes (\Sp(n)\cdot \Sp(1)\times \RR^+)=
\Aut_{qc}(\cM)$ up to conjugate.
In particular the closed subgroup $G$ is also amenable so that
$G$ is an extension of a solvable group by a compact group.
Let $G=R\rtimes K$ be the semidirect product 
where $R$ is solvable and $K$ is compact.
It is noted from \cite{OK2} that there is a simply connected characteristic
solvable subgroup $R_0$ such that $R=R_0\cdot T$ where
$T$ is a maximal compact subgroup of $R$. Letting $H=T\cdot K$,
we have a semidirect product
$G=R_0\rtimes H$. 

Let $\Psh_{qc}(\cM)=\cM\rtimes (\Sp(n)\cdot \Sp(1))\leq \Aut_{qc}(\cM)$ and
$\cM=S^{4n+3}-\{\infty\}$ be as above.
Suppose $\displaystyle M=\cG/\cK$ is a simply connected homogeneous $qc$-manifold
where $\cG\leq \Aut_{qc}(M)$ and $\cK$ is a closed subgroup of $\cG$.
Let $\displaystyle (\rho, \dev): (\Aut_{qc}(M),M)\ra (\Aut_{qc}(\cM),\cM))$
be the developing pair. Put $\rho(\cG)=G$.

\begin{pro}\label{tran3}
If $G\leq \Psh_{qc}(\cM)$ up to conjugate, then
$\dev(M)=\cM$ for which $M$ is a homogeneous 
$qc$-manifold by the $qc$-group $\cM\rtimes (\Sp(n)\cdot \Sp(1))$.
\end{pro}

\begin{proof}
The developing pair reduces to
$\displaystyle (\rho, \dev): (\cG,M)\ra (\Psh_{qc}(\cM),\cM)$.
Taking an $\cM\rtimes (\Sp(n)\cdot \Sp(1))$-invariant Riemannian metric 
on $\cM$, the pullback metric by $\dev$
induces a homogeneous Riemannian metric on $M$.
Noting $M$ is geodesically complete (\cf \cite{WO}),
$\dev:M\ra \cM$ is an isometry. Hence $M$ is $qc$-homogeneous by the group 
$\cM\rtimes (\Sp(n)\cdot \Sp(1))$. 
\end{proof}

In general case for $G=R_0\rtimes H\leq \Aut_{qc}(\cM)$,
consider the exact sequence:
\begin{equation}\label{trancom}
\begin{CD}
1\ \ra\,\cM\rtimes(\Sp(n)\cdot \Sp(1))@>>>\Aut_{qc}(\cM)@>{\rm pr}>>\RR^+\ \ra\ 1\\
\end{CD}
\end{equation}where ${\rm pr}(G)=\RR^+$ from Proposition \ref{tran3}.
Then there exists a one-parameter subgroup $A\leq R_0$
such that ${\rm pr}(A)=\RR^+$. It follows
$R_0=N\rtimes A$
where $N$ is a nilpotent subgroup
such that $\displaystyle N=(\ker\, {\rm pr})\cap R_0=\cM\cap R_0$.
Then $G=(N\rtimes A)\cdot H$.
Since $A\leq \PU(n+1,1)$ is of non-elliptic type which
stabilizes $\{\mbi{0},\infty\}\subset S^{4n+3}$ up to conjugate,
there is a geodesic segment between 
$\mbi{0}$ and $\infty$ translated by $A$. As $H\leq \Sp(n)\cdot \Sp(1)$ 
fixes $\{\mbi{0},\infty\}$, 
the conjugate $h\cdot A\cdot h^{-1}$ for $h\in H$
fixes $\{\mbi{0},\infty\}$ also, there is a geodesic segment between 
$\mbi{0}$ and $\infty$ translated by $h\cdot A\cdot h^{-1}$. Then
these two geodesic segments between the same endpoints spans a
(geodesically) flat plane in $\HH^{n+1}_\HH$.
Since $\HH^{n+1}_\HH$ is a complete simply connected
Riemannian manifold of constant negative quaternionic curvature,
there is no such flat plane so that $h\cdot A\cdot h^{-1}=A$.
In particular, the elements of $A$ and $H$ commute.
We have 
\begin{equation}\label{radicalA}
G=N\rtimes (A\times H).
\end{equation}
Let $q:\cM\ra \HH^n$ be the projection which is equivariant
with respect to the homomorphism ${\sf q}$ in the exact sequence:
\begin{equation*}
\begin{CD}
1\ \ra\,\RR^3\lra \cM\rtimes(\Sp(n)\cdot \Sp(1)\times \RR^+)\stackrel{{\sf q}}\lra\HH^n\rtimes(\Sp(n)\cdot \Sp(1)\times \RR^+)\ \ra\ 1\\
\end{CD}
\end{equation*} where ${\sf q}(G)={\sf q}(N)\rtimes (A\times H)$ by \eqref{radicalA}. 
Denote the point $\mbi{1}$ in $\cM=\RR^3\times \HH^n$ by 
\begin{equation*}\begin{split}
\mbi{1}=\left[\mbi{0},\binom 1{\mbi{0}} \right]\, \mbox{where} \ 
\binom 1{\mbi{0}}
={}^t(\,\overbrace{\,1,0,\cdots,0,\,}^{k} \,\overbrace{\,0,\ldots,0\,}^{n-k}\,).
\end{split}\end{equation*}
Thus $\displaystyle q(\mbi{1})=
\binom 1{\mbi{0}}\in \HH^n$.
Since $\sfD_0$ on $\cM$ maps to $q_*(\sfD_0)=T\HH^n$,
the homogeneous quaternionic manifold ${\sf q}(G)\cdot \mbi{0}$ 
becomes $\displaystyle \binom{(A\times H)\cdot 1}{{\sf q}(N)\cdot 0}\subset \HH^n$.
Here either $\displaystyle (A\times H)\cdot 1=
A\times H/H_{1}=\RR^+\times S^{4k-1}=\HH^{k}-\{\mbi{0}\}$
 and $\displaystyle {\sf q}(N)\cdot 0=\HH^{n-k}$ $(1\leq k\leq n)$, or
$\displaystyle (A\times H)\cdot{1}=A\cdot{1}=\RR^+$ and
$\displaystyle {\sf q}(N)\cdot 0=\binom {{\rm Im} \HH}{\HH^{n-1}}$.

\begin{lemma}\label{summary}
The orbit $G\cdot \mbi{1}$ at $\mbi{1}\in \cM$ becomes: 
\[\RR^3\times \left[\begin{array}{c}
\HH^{k}-\{\mbi{0}\}\\
  \HH^{n-k}\end{array}\right] 
 \   (1\leq k\leq n) \ \  \mbox{or}\ \
\RR^3\times \left[\begin{array}{c}
{\rm Im}\HH +\RR^+\\
\HH^{n-1}\\
\end{array}\right].
\]
\end{lemma}

\begin{proof}
According to each orbit, $G$ is isomorphic to either $G(k)$ or $G(0)$
respectively:
\begin{enumerate}
\item $G(k)=\left( \RR^3\cdot\left[\begin{array}{c}
\mbi{0}\\
\HH^{n-k}
 \end{array}\right]\right)\rtimes 
\left(\begin{array}{cc}
\Sp(k)& \\
 &\Sp(n-k)
\end{array}\right)\cdot \Sp(1)\times \RR^+$ \\  $(1\leq k\leq n)$.

\item $G(0)=\left( \RR^3\cdot\left[\begin{array}{c}
{\rm Im}\,\HH\\
\HH^{n-1}\\
\end{array}\right]\right)\rtimes 
\left(\begin{array}{cc}
a& \\ 
&\Sp(n-1) \\
\end{array}\right)\cdot a\times \RR^+$\  $(a\in \Sp(1))$.
\end{enumerate}
For this, an element of $G(k)$ acts at $\mbi{1}$ as 
\begin{equation*}\begin{split}
&(\mbi{t},\left[\begin{array}{c}
\mbi{0}\\
z\\
\end{array}\right])\cdot 
\left(\left(\begin{array}{cc}
A& \\
 &B \end{array}\right)\cdot \lam\cdot a \right)
\cdot \left[\mbi{0},\binom 1{\mbi{0}} \right]
=(\mbi{t},\left[\begin{array}{c}
\mbi{0}\\
z\\
\end{array}\right])\cdot (\mbi{0},\left[\begin{array}{c}
\lam A\cdot 1\cdot \bar a\\
\mbi{0}\\
\end{array}\right])\\
&\quad \hskip6.4cm \quad =
(\mbi{t},\left[\begin{array}{c}
\lam A\cdot 1\cdot \bar a\\
z\\
\end{array}\right]).\end{split}
\end{equation*}

Noting $\lam A\cdot 1\cdot \bar a\neq \mbi{0}$,
it is easy to see 
\begin{equation}\label{Gk}
G(k)\cdot \mbi{1}=\RR^3\times \left[\begin{array}{c}
\HH^k-\{\mbi{0}\}\\
\HH^{n-k}\\
\end{array}\right].
\end{equation}

Similarly the orbit of $G(0)$ at $\mbi{1}$ consists of
\begin{equation*}\begin{split}
(\mbi{t},\left[\begin{array}{c}
u\\
z\\
\end{array}\right])\cdot 
\left(\left(\begin{array}{cc}
a& \\
 &A \end{array}\right)\cdot \lam\cdot a \right)
\cdot &\left[\mbi{0},\binom 1{\mbi{0}} \right]=
(\mbi{t},\left[\begin{array}{c}
u\\
z\\
\end{array}\right])\cdot 
\left[\mbi{0},\binom {\lam a\cdot 1\cdot \bar a}{\mbi{0}} \right]\\
&=(\mbi{t},\left[\begin{array}{c}
u+\lam a\cdot 1\cdot \bar a\\
z\\
\end{array}\right])
=
(\mbi{t},\left[\begin{array}{c}
u+\lam\\
z\\
\end{array}\right])
\end{split}
\end{equation*}
such that ${\rm Re}(u+\lam)=\lam>0$.
Putting $\HH_+=\{z\in \HH\mid {\rm Re}(z)>0\}$, it follows
\begin{equation}\label{G0}
G(0)\cdot \mbi{1}=\RR^3\times \left[\begin{array}{c}
\HH_+\\
\HH^{n-1}\\
\end{array}\right].
\end{equation}
\end{proof}

\noindent Let $\displaystyle G(k)_{\mbi{1}}$
be the stabilizer at $\mbi{1}$. Suppose
$\displaystyle (\mbi{t},\left[\begin{array}{c}
\mbi{0}\\
z\\
\end{array}\right])\cdot 
\left(\left(\begin{array}{cc}
A& \\
 &B \end{array}\right)\cdot \lam\cdot a \right)
\cdot \left[\mbi{0},\binom 1{\mbi{0}} \right]=
\left[\mbi{0},\binom 1{\mbi{0}} \right]$.
Then
$\mbi{t}=\mbi{0}, z=\mbi{0}$. 
The equation
$\lam A\cdot 1\cdot \bar a=1$ shows 
$\lam |A\cdot 1\cdot \bar a|=1$ where $\displaystyle 1={}^t[\,\overbrace{\,1,0,\ldots,0\,}^{k}\,]$. So $\lam =1$.
As $\displaystyle A\cdot 1\cdot \bar a=
{}^t[a_{11}\cdot \bar a,a_{21}\cdot \bar a,\ldots,a_{k1}\cdot \bar a]=
{}^t[\,1,0,\ldots,0\,]$, it follows $a_{11}=a,\, a_{21}=\cdots=a_{k1}=0$.
Then
$\displaystyle A\cdot a=\left(\begin{array}{cc}
a & 0\\
0 & A_{k-1} \\
\end{array}\right)\cdot a \ \ ({}^\exists\, A_{k-1}\in \Sp(k-1))$.
Letting $a \in\Sp(1)$,

\begin{equation}\label{gk1comp}
G(k)_{\mbi{1}}=\left(\mbi{0},\left[\begin{array}{c}
\mbi{0}\\
\mbi{0} \end{array}\right]\right)\rtimes 
\left(\begin{array}{ccc}
a& &\\
& \Sp(k-1) & \\
& &\Sp(n-k) \\
\end{array}\right)\cdot a
\end{equation}
Similarly if $g\cdot \mbi{1}=\mbi{1}$ for $g\in G(0)$, then 
note ${\rm Re}(u+\lam)=\lam=1$ from \eqref{G0}. As $u+1=1$, $u=0$.
It follows
\begin{equation}\label{g01comp}
G(0)_{\mbi{1}}=\left(\mbi{0},\left[\begin{array}{c}
\mbi{0}\\
\mbi{0} \end{array}\right]\right)\rtimes 
\left(\begin{array}{cc}
 a & \\
&\Sp(n-1) \\
\end{array}\right)\cdot a\ \ (a\in \Sp(1)).
\end{equation}
Since both $G(k)_{\mbi{1}}$ and $G(0)_{\mbi{1}}$ are compact subgroups,
we obtain 
\begin{pro}\label{homGk}
The $qc$-manifolds $\displaystyle G(k)/G(k)_{\mbi{1}},\ G(0)/G(0)_{\mbi{1}}$
are homogeneous Riemannian domains in $\cM$: 
\begin{enumerate}
\item $\displaystyle G(k)/G(k)_{\mbi{1}}=G(k)\cdot{\mbi{1}}=\RR^3\times \left[\begin{array}{c}
\HH^{k}-\{\mbi{0}\}\\
  \HH^{n-k}\end{array}\right]
 \\  (1\leq k\leq n)$. 
\item $\displaystyle G(0)/G(0)_{\mbi{1}}=G(0)\cdot{\mbi{1}}=\RR^3\times \left[\begin{array}{c}
\HH_+\\ 
\HH^{n-1}\\
\end{array}\right]$.
\end{enumerate}
\end{pro}
Let $X(k)$ $(0\leq k\leq n)$ be
the orbit $G(k)\cdot\mbi{1}$ $(k\neq 0)$
or $G(0)\cdot\mbi{1}$. 
We obtain the following. (Compare Proposition \ref{homGk}, Proposition \ref{tran2}, \eqref{R3bundle}.) 

\begin{theorem}\label{isoqchomo}
Any simply connected spherical homogeneous $qc$-manifold $M$ is $qc$-isomorphic to
$S^{4n+3}$, $S^{4n+3}-S^{4m-1}$ $(1\leq m\leq n)$, 
or $\cM$, $X(k)$ $(0\leq k\leq n)$. 
In particular, only $S^{4n+3}$ admits the $qck$-group $\Sp(1)$ 
and each $\cM$ or $X(k)$
admits the $qck$-group $\RR^3$.
\end{theorem}
\begin{proof}The proof divides into two cases whether the holonomy image
$G=\rho(\cG)\leq \PSp(n+1,1)$ is semisimple or not. Noting
$M$ is simply connected, there is a $\rho$-equivariant
developing map $\displaystyle \dev:M\ra S^{4n+3}$. 
Put $X=\dev(M)=G\cdot p$.
When $G$ is semisimple, Proposition \ref{tran2} shows
the only homogeneous $qc$-manifold $X$ is
$\displaystyle S^{4n+3}-S^{4m-1}=
{\rm Sp}(m,1)\cdot{\rm Sp}(n-m+1)/
({\rm Sp}(m)\times (\Delta\Sp(1)\cdot {\rm Sp}(n-m))$\ 
$(1\leq m\leq n)$ (\cf \cite[Lemma 3.3 (ii)]{KU}).
In addition only $X=S^{4n+3}=\Sp(n+1)/\Sp(n)$ admits the $qck$-group $\cR=\Sp(1)$.
Since all of these $X$ are Riemannian homogeneous, 
the pullback metric by $\dev$ implies $M$ is $qc$-isomorphic to $X$ as before.
When $G$ has the nontrivial radical, $G$ is 
either $G(k)$ or $G(0)$ by Proposition \ref{homGk}.
Let $A$ be the one-parameter subgroup of $G$ from \eqref{radicalA}.
Since $A$ stabilizes $\{\mbi{0},\infty\}$ in $S^{4n+3}$,
the limit set $L(G)$ contains $\{\mbi{0},\infty\}$ (\cf \cite{CG}). 
The complement
$S^{4n+3}-X$ is a closed subset containing more than one point
(otherwise $X= S^{4n+3}-\{\infty\}=\cM)$.
Then it follows from \cite{CG} that
$L(G)\subset S^{4n+3}-X$, that is 
\begin{equation}\label{limimal}
X\subset S^{4n+3}-L(G).
\end{equation}
Let $\displaystyle \mbi{0}=\left[\mbi{0},\binom 00 \right]$ be the origin of
$\cM$. According to whether 
$G$ is $G(k)$ or $G(0)$,
the orbit $G\cdot \mbi{0}$ becomes
(i)\, 
$\displaystyle G(k)\cdot \mbi{0}
=
 \RR^3\times\left[\begin{array}{c}
\mbi{0}\\
\HH^{n-k}
 \end{array}\right]$,
 \ (ii)\,\  
$\displaystyle G(0)\cdot \mbi{0}
=\RR^3\times \left[\begin{array}{c}
{\rm Im}\,\HH\\
\HH^{n-1}\\
\end{array}\right]$.
In each case the union $G\cdot\mbi{0}\,\cup\{\infty\}$
is a closed subset in $S^{4n+3}$. (In fact, it is the
sphere \emph{diffeomorphic} to either
$\displaystyle G(k)\cdot \mbi{0}\,\cup\{\infty\}=S^{4(n-k)+3}$ or
$\displaystyle G(0)\cdot\mbi{0}\,\cup\{\infty\}=S^{4n+2}$.)
Since the limit set $L(G)$ is $G$-invariant
with $\{\mbi{0},\infty\}\subset L(G)$, it follows
$L(G)\subset G\cdot \mbi{0}\,\cup
\{\infty\}\subset L(G)$, that is  $\displaystyle L(G)=G\cdot \mbi{0}\,\cup \{\infty\}$.
From \eqref{limimal},
\begin{equation}\label{limieq}
X\subset S^{4n+3}-L(G)=\cM\cup\{\infty\}-(G\cdot \mbi{0}\cup \{\infty\})=
\cM-G\cdot \mbi{0}.
\end{equation}
Calculate for $G=G(k)$,
\begin{equation*}\begin{split}
\cM-G(k)\cdot \mbi{0}&=\RR^3\times \HH^n- \RR^3\times\left[\begin{array}{c}
\mbi{0}\\
\HH^{n-k}
 \end{array}\right]
=\RR^3\times
\left[\begin{array}{c}
\HH^k-\{\mbi{0}\}\\
\HH^{n-k}
 \end{array}\right].\\
\cM-G(0)\cdot \mbi{0}&=\RR^3\times \HH^n- \RR^3\times\left[\begin{array}{c}
{\rm Im}\,\HH\\
\HH^{n-1}\\
\end{array}\right]
=\RR^3\times 
\left[\begin{array}{c}
\HH_+\cup \HH_-\\
\HH^{n-1}
 \end{array}\right].
\end{split}\end{equation*}
It follows from \eqref{Gk} and \eqref{G0} that
\begin{equation}
\begin{split}
 X\subset S^{4n+3}-L(G(k))&=G(k)\cdot \mbi{1},\\
 X\subset (S^{4n+3}-L(G(0)))^0&=G(0)\cdot \mbi{1}.\\
\end{split}
\end{equation}
By Proposition \ref{homGk}, $G\cdot \mbi{1}$ is homogeneous Riemannian.
Since $\dev(M)=X$ is homogeneous, $X=G\cdot \mbi{1}$.
As above $M$ is $qc$-isometric to $X=G\cdot\mbi{1}$.
\end{proof}

By Theorem \ref{isoqchomo} (\cf Proposition \ref{homGk}), it has a principal bundle:
\begin{equation}\label{R3bundle}\begin{CD}
\RR^3@>>>X(k)@>>>Y(k)=\left[\begin{array}{c}
\HH^{k}-\{\mbi{0}\}\\
  \HH^{n-k}\
\end{array}\right],\\
\RR^3@>>> X(0)@>>>Y(0)=\left[\begin{array}{c}
\HH_+\\
\HH^{n-1}\\
\end{array}\right]
\end{CD}\end{equation}where $Y(k)$ is a domain of $\HH^n$.   
The standard $qc$-form $\om_0$ on $\cM$
restricts an invariant $qc$-structure to $X(k)$.
Since $G(k)$ has $\langle \lam\rangle=\RR^+$,   
$\lam^*\om_0=\lam^2\om_0$ by \eqref{actionqc},
though $G(k)$ preserves the $qc$-structure.
Thus $\Aut_{qc}(X(k))=G(k)$. 
Putting $\Psh(X(k))=\Psh_{qc}(X(k),\om_0,\{J_\al\}_{\al=1}^3)$
$(0\leq k\leq n)$, 

\begin{equation*}\begin{split}
\Psh_{qc}(X(k))&=
\left( \RR^3\cdot\left[\begin{array}{c}
\mbi{0}\\
\HH^{n-k}
 \end{array}\right]\right)\rtimes 
\left(\begin{array}{cc}
\Sp(k)& \\
 &\Sp(n-k)
\end{array}\right)\cdot \Sp(1),\\ 
\Psh_{qc}(X(0))&=
\left( \RR^3\cdot\left[\begin{array}{c}
{\rm Im}\,\HH\\
\HH^{n-1}\\
\end{array}\right]\right)\rtimes 
\left(\begin{array}{cc}
a& \\ 
&\Sp(n-1) \\
\end{array}\right)\cdot a \ \ ({}^\forall\,a\in \Sp(1)).\\ 
\end{split}
\end{equation*}However
$\displaystyle \Psh_{qc}(X(k),\om_0,\{J_\al\}_{\al=1}^3)$ is not transitive on $X(k)$.
Of course, there is a $qc$-form $\eta$ 
such that $\displaystyle \Psh_{qc}(X(k),\eta,\{J_\al\}_{\al=1}^3)=
\Aut_{qc}(X(k))=G(k)$ by Theorem \ref{qvanish}. 
Since $\eta=v\cdot \om_0$ for some non-constant $v\in C^\infty(X(k),\RR^+)$,
$\RR^3$ does not induce a $qck$-distribution for $\eta$.

\section{Curvature criterion of $qc$-manifolds with $qck$-group $\RR^3$}\label{app}
Let $\displaystyle (X,g,\eta,\{J_\al\}_{\al=1}^{3})$ be
a $4n+3$-dimensional simply connected non-compact positive definite 
$qc$-manifold with $qck$-distribution
$V=\RR^3$ where $g=g_\eta$ (\cf \eqref{Riemann}). $(X,g)$ is
$\sfD$-Einstein, that is

\begin{pro}\label{Drist}
$\displaystyle {\rm Ric}(\mbi{x},\mbi{y})=-6g(\mbi{x},\mbi{y})$\
$({}^\forall\,\mbi{x},\mbi{y}\in \sfD)$.
\end{pro}

\begin{proof}
Let $\displaystyle\sfD=\{\mbi{v}_1,\ldots,\mbi{v}_{4n}\}=
\{\mbi{v}_1,\ldots,\mbi{v}_n, J_\al\mbi{v}_1,\ldots,
J_\al\mbi{v}_n,\ \al=1,2,3\}$. Choose $\mbi{x},\mbi{y}\in\sfD$ satisfying
$g(\mbi{x},\mbi{x})=g(\mbi{y},\mbi{y})=1$.
O'Neill's formula \cite[(3.30)]{CE} shows ($\mbi{x},\mbi{y},\mbi{v}\in \sfD$)\,:
\begin{equation}\label{ONEL}
g(R(\mbi{v},\mbi{x})\mbi{y},\mbi{v})=
\hat g(\hat R(\hat{\mbi{v}},\hat{\mbi{x}})\hat{\mbi{y}},\hat{\mbi{v}})+
\frac 34g([\mbi{v},\mbi{y}]^V,[\mbi{x},\mbi{v}]^{V}).
\end{equation} 
As $g=g_\eta$ in \eqref{Riemann},
\begin{equation*}\begin{split}
&\sum_{i=1}^{4n}\frac 34\bigl(\sum_{\al=1}^3\eta_\al([\mbi{v}_i,\mbi{y}]^V)
\eta_\al([\mbi{x},\mbi{v}_i]^V)\bigr)
=-3\sum_{i=1}^{4n}\bigl(\sum_{\al=1}^3 d\eta_\al(\mbi{v}_i,\mbi{y})d\eta_\al(\mbi{v}_i,\mbi{x})\bigr)\\
&= -3\sum_{\al=1}^3\,d\eta_\al(J_\al\mbi{x},\mbi{y})d\eta_\al(J_\al\mbi{x},\mbi{x})
=-3\cdot 3g(\mbi{x},\mbi{y})=-9g(\mbi{x},\mbi{y}).
\end{split}\end{equation*}
Consider the $qc$-bundle \eqref{BW3} (Riemannian submersion):
$\displaystyle \RR^{3}\ra\,(X,g)\stackrel{p}\lra\, (Y,\hat g)$
where $\langle \xi_1, \xi_2,\xi_3\rangle\oplus {\sfD}=TX$.
 As $Y$ is hyperK\"ahler,
$\displaystyle 
{\rm Ric}(\hat{\mbi{x}},\hat{\mbi{y}})=\sum_{i=1}^{4n}\hat g\bigl(\hat R(\hat{\mbi{v}}_i,\hat{\mbi{x}})\hat{\mbi{y}},\hat{\mbi{v}}_i\bigr)=0$.
 Noting \eqref{ONEL}, 
\begin{equation}\label{4nqc}
\sum_{i=1}^{4n}g(R(\mbi{v}_i,\mbi{x})\mbi{y},\mbi{v}_i)=-9g(\mbi{x},\mbi{y}).
\end{equation}

On the other hand, 
take $\RR^2=\{\xi_2,\xi_3\}$ for which 
there is the Riemannian submersion
$\displaystyle \RR^{2}\ra\,(X,g)\stackrel{\mu}\lra (X_2,g_2)$
where
$(X_2,\hat\eta_1,J_1',\hat\xi_1)$ with $\mu_*\xi_1=\hat\xi_1$
is the pseudo-Hermitian (Sasaki) manifold such that
\begin{equation}\label{g2met}
g_2(\hat {\mbi u},\hat{\mbi v})=
\hat\eta_1(\hat {\mbi u})\cdot\hat\eta_1(\hat {\mbi v})+
d\hat\eta_1(J_1'\hat {\mbi u},\hat {\mbi v})\,
\ (\hat {\mbi u},\hat{\mbi v}\in TX_2). 
\end{equation}
Put ${\sf E}=\langle\xi_1\rangle\oplus \sfD\subset TX$. Then
$\displaystyle \mu_*: {\sf E}\lra\,TX_2=\langle\hat\xi_1\rangle\oplus\mu_*\sfD$
is an isometry such that $\langle\xi_2,\xi_3\rangle={\sf E}^\perp$. 
Let $R_2$ be the Riemannian curvature on $(X_2,g_2)$. 
In general, a Sasaki manifold $(X_2,g_2,\hat\eta_1,J_1')$ satisfies
$\displaystyle R_2(\hat \xi_1,\hat{\mbi{x}})\hat{\mbi{y}}
=g_2(\hat{\mbi{x}},\hat{\mbi{y}})\hat\xi_1-\hat\eta_1(\hat{\mbi{y}})\hat{\mbi{x}}\
\, (\hat{\mbi{x}},\hat{\mbi{y}}\in TX_2)$. (See \cite[(2.5)]{ST} for instance.) 
If $\mbi{x},\mbi{y}\in \sfD\subset TX$ with $\mu_*(\mbi{x})=\hat{\mbi{x}}$,
$\mu_*(\mbi{y})=\hat{\mbi{y}}$, then 
$g(\mbi{x},\mbi{y})=d\eta_1(J_1\mbi{x},\mbi{y})=g_2(\hat{\mbi{x}},\hat{\mbi{y}})$ by \eqref{g2met} (\cf \eqref{hKaehler}) and so the above equation becomes
$\displaystyle R_2(\hat \xi_1,\hat{\mbi{x}})\hat{\mbi{y}}=g(\mbi{x},\mbi{y})\hat\xi_1$.
O'Neill's formula shows 
$\displaystyle g(R(\xi_1,\mbi{x})\mbi{y},\xi_1)=
g_2(R_2(\hat\xi_1,\hat{\mbi{x}})\hat{\mbi{y}},\hat\xi_1)$
since $\displaystyle [\xi_1,\mbi{x}]^{\cV}=0$ where $\cV=\{\xi_2,\xi_3\}$.
Noting $g_2(\hat\xi_1,\hat\xi_1)=1$, 
substitute the above equation:\begin{equation}\label{cal1}
g(R(\xi_1,\mbi{x})\mbi{y},\xi_1)=g(\mbi{x},\mbi{y})\ \, (\mbi{x},\mbi{y}\in \sfD).
\end{equation}
Applying this argument to the other submersion
$\displaystyle \RR^{2}\ra (X,g) \stackrel{\mu'}\lra (X_2',g_2')$
for which $(X'_2,\hat\eta_2,J_2',\hat \xi_2)$ 
is the pseudo-Hermitian (Sasaki) manifold with $\RR^{2}=\{\xi_3,\xi_1\}$
\,(respectively so is
$(X''_2,\hat\eta_3,J_3',\hat \xi_3)$\, with $\RR^{2}=\{\xi_1,\xi_2\}$\,).
Similarly 
\begin{equation}\label{cal2}
g(R(\xi_2,\mbi{x})\mbi{y},\xi_2)=g(\mbi{x},\mbi{y}),\ \
g(R(\xi_3,\mbi{x})\mbi{y},\xi_3)=g(\mbi{x},\mbi{y})\ (\mbi{x},\mbi{y}\in \sfD).
\end{equation}
Using \eqref{4nqc}, \eqref{cal1}, \eqref{cal2},
for any $\mbi{x},\mbi{y}\in \sfD$,
 we obtain
\begin{equation}\label{ricsum1}\begin{split}
{\rm Ric}(\mbi{x},\mbi{y})&=\sum_{i=1}^3 g(R(\xi_i,\mbi{x})\mbi{y},\xi_i)+
\sum_{i=1}^{4n}g(R(\mbi{v}_i,\mbi{x})\mbi{y},\mbi{v}_i)=-6g(\mbi{x},\mbi{y}).
\end{split}\end{equation}
\end{proof}

\begin{pro}\label{notEinstein}
Let $(X,g,\eta,\{J_\al\}_{\al=1}^3)$ be as above.
Unlike the $3$-Sasaki structure, $(X,g)$ is not Einstein. 
Indeed we obtain 
\begin{equation}\label{ricxi}
{\rm Ric}\,(\xi_\al,\xi_\al)=4n\ \ (\al=1,2,3).
\end{equation}
\end{pro}

\begin{proof}
Let $\displaystyle \RR^{2}=\{\xi_2,\xi_3\}\ra\, (X,g)\stackrel{\mu}\lra (X_2,g_2)$
be the Riemannian submersion as in \eqref{g2met}.
Take an orthonormal basis 
$\{\mbi{u}_1,\ldots,\mbi{u}_{4n+1}\}=\mbi{W}$
such that 
$\displaystyle\langle\xi_2,\xi_3\rangle\oplus \mbi{W}\,
(=\langle\xi_1,\xi_2,\xi_3\rangle\oplus \sfD)$
forms a basis of $TX$. Then
${\rm Ric}\, (\xi_1,\mbi{x})=\sum_{i=1}^{4n+1}g(R(\mbi{u}_i,\xi_1)\mbi{x},\mbi{u}_i)+
\sum_{\al=2}^{3}g(R(\xi_\al,\xi_1)\mbi{x},\xi_\al)$ $(\xi_1,\mbi{x}\in TX)$.
Since $\xi_1$ is Killing, 
$[\xi_1,\mbi{u}_i]^{\cV}=\mbi{0}$ where $\cV=\langle\xi_2,\xi_3\rangle$.
O'Neill's formula implies
\begin{equation}\label{g2reduc}
g(R(\mbi{u}_i,\xi_1)\mbi{x},\mbi{u}_i)=
g_2(R_2(\hat{\mbi{u}}_i,\hat\xi_1)\hat{\mbi{x}},\hat{\mbi{u}}_i).
\end{equation}
As $(X_2,g_2, \hat\eta_1, J_1',\hat\xi_1)$
is a $4n+1$-dimensional Sasaki manifold, 
the fundamental property of Sasaki manifold 
shows
$\displaystyle {\rm Ric}\,(\hat\xi_1,\hat{\mbi{x}})=4n g_2(\hat\xi_1,\hat x)$,
 (see \cite[(1.6)]{ST2} for example.)
 By \eqref{g2reduc}, 
$\displaystyle {\rm Ric}\,(\xi_1,\mbi{x})
=4n g_2(\hat{\xi}_1,\hat{\mbi{x}})+g(R(\xi_2,\xi_1)\mbi{x},\xi_2)+
g(R(\xi_3,\xi_1)\mbi{x},\xi_3)$.
Replace $\mbi{x}$ by $\xi_1$.
Noting $\RR^3=\langle \xi_1,\xi_2,\xi_3\rangle$ spans a flat geodesic subspace of $X$,
$g(R(\xi_2,\xi_1)\xi_1,\xi_2)=K(\xi_2,\xi_1)=0$,
$g(R(\xi_3,\xi_1)\xi_1,\xi_3)=K(\xi_3,\xi_1)=0$.
Hence
$\displaystyle  {\rm Ric}\,(\xi_1,\xi_1)=4n g_2(\hat\xi_1,\hat\xi_1)=4n$.
Apply the same argument to
each Sasaki manifold
 $\displaystyle (X_2',g_2',\hat\eta_2,J_2',\hat\xi_2)$,
$\displaystyle (X_2'',g_2'',\hat\eta_3,J_3',\hat\xi_3)$, it follows
$\displaystyle {\rm Ric}\,(\xi_2,\xi_2)=4n g(\xi_2,\xi_2)=4n$,
$\displaystyle {\rm Ric}\,(\xi_3,\xi_3)=4n g(\xi_3,\xi_3)=4n$ respectively.
\end{proof}

\begin{remark}
By Theorem $\ref{hcr}$, letting $X_1=X/\RR, X_2=X/\RR^2$,
there is a \emph{holomorphic bundle} 
$\displaystyle \RR^{2}\ra (X_1,g_1)\stackrel{q}\lra (Y,\hat g)$, 
a \emph{pseudo-Hermitian Sasaki bundle}
$\displaystyle \RR\ra (X_2,g_2)\stackrel{\nu}\lra (Y,\hat g)$
respectively. The similar argument shows that
${\bf (i)}$ $(X_1,g_1)$ is $p_{1*}{\sfD}$-Einstein, that is 
$\displaystyle {\rm Ric}(\hat{\mbi{x}},\hat{\mbi{y}})=
-4g_1(\hat{\mbi{x}},\hat{\mbi{y}})$\, $(\hat{\mbi{x}},\hat{\mbi{y}}\in p_{1*}\sfD)$
where  $\langle \xi_2,\xi_3\rangle\oplus p_{1*}{\sfD}=TX_1$,
${\bf (ii)}$ $(X_2,g_2)$ is $\mu_*{\sfD}$-Einstein,
$\displaystyle {\rm Ric}(\hat{\mbi{x}},\hat{\mbi{y}})
=-2g_2(\hat{\mbi{x}},\hat{\mbi{y}})$\ \,  $(\hat{\mbi{x}},\hat{\mbi{y}}\in \mu_*\sfD)$
where $\langle \hat\xi_1\rangle\oplus\mu_*\sfD=TX_2$.
\end{remark}

\begin{pro}\label{Dristen}
Let $\displaystyle (X,g,\eta,\{J_\al\}_{\al=1}^{3})$ be
a $4n+3\,(\geq 11)$-dimensional
simply connected non-compact positive definite $qc$-manifold with $qck$-group
$\RR^3$. 
For some $\al\in \{1,2,3\}$,
if the sectional curvature $K(\mbi{x},J_\al\mbi{x})$ is \emph{constant} for every unit vector $\mbi{x}\in \sfD$ at all points of $X$,
then there is a $qc$-isometric immersion $\dev
: (X,g)\ra\, (\cM,g_0)$ with $\dev^*g_0=g$
such that $\displaystyle \dev^*\om_0=\eta$, $\dev_*J_\al=J_\al\dev_*|_{\sfD}$ and 
$\displaystyle\dev \circ\, t=t\circ \dev$ $({}^\forall\, t\in \RR^3)$.
Moreover, if there exists a discrete subgroup $\Gamma\leq \Psh_{qc}(X,\eta)$
such that $M=X/\Gamma$ is compact, then $\dev$ induces
a $qc$-isometry of $M$ onto the quaternionic Heisenberg infranilmanifold $\cM/\rho(\Gamma)$ 
for the holonomy group $\rho(\Gamma)\leq \cM\rtimes(\Sp(n)\cdot \Sp(1))$.
\end{pro}

\begin{proof}
Let $\displaystyle \RR^{3}\ra (X,g)\stackrel{p}\lra (Y,\hat g)$
be a Riemannian submersion where $(Y,\hat g)$ is a hyperK\"ahler manifold as before. 
In particular, each induced endomorphism $\hat J_\al$ is a complex structure on $Y$.
For every unit vector $\displaystyle \hat{\mbi{x}}\in TY$,
choose $\mbi{x}\in \sfD$ such that $\displaystyle p_*\mbi{x}=\hat{\mbi{x}}$. 
By \cite[(3.20) Theorem]{CE}, the \emph{holomorphic} sectional curvature 
has the relation $\displaystyle  K(\mbi{x},J_\al\mbi{x})
=\hat K(\hat{\mbi{x}},\hat J_\al\hat{\mbi{x}})-3$.
As $K(\mbi{x},J_\al\mbi{x})$ is constant by the hypothesis,
so is $\displaystyle \hat K(\hat{\mbi{x}},\hat J_\al\hat{\mbi{x}})$ on $Y$.
Let $\displaystyle \hat K(\hat{\mbi{x}},\hat J_\al\hat{\mbi{x}})=c$
for every unit vector $\mbi{x}\in \sfD$ at all points of $Y$.
Then $\displaystyle (Y,\hat g, \hat J_\al)$ is 
locally holomorphically isometric to
a K\"ahler complex space form of constant curvature $c$. Since
$\displaystyle (Y,\hat g,\{\hat J_\al\}_{\al=1}^3)$ is hyperK\"ahler,
$c=0$. It follows from \cite[II. Proposition 7.3\, IX]{KN}
that $(Y,g)$ is of zero curvature and hence 
$Y$ is locally isometric to $\HH^n$ for ${\rm dim}\, Y\geq 8$ by the uniformization 
(\cf \cite[Theorem 3.9, also Theorem 5.2, Theorem 3.5]{IS}, \cite[Corollary 3]{AL}).
Let $\displaystyle \Dev :Y\ra \HH^n$ be a quaternionically 
isometric immersion.
If $\displaystyle (g_0,\Om_0,\{{\hat J}_\al\}_{\al=1}^3)$ is
the standard quaternionic structure, then
$\displaystyle \Dev$ satisfies
$\displaystyle \Dev^*\Om_0=\Om$ and 
$\displaystyle\Dev_*\circ{\hat J}_\al={\hat J}_\al\circ\Dev_*$
up to conjugate by an element of $\Sp(1)$ $(\al=1,2,3)$. Since 
$\displaystyle \RR^3\ra X\lra \, Y$ is the principal $qc$-bundle
as well as the standard principal $qc$-bundle
$\displaystyle \RR^3=C(\cM)\ra \cM\lra\, \HH^n$,
applying Proposition \ref{exact3} to $\Dev$,
there is a lift of $qc$-immersion $\dev:X\ra \cM$ 
satisfying $\dev^*\om_0=\eta,\, \dev_*\circ{J}_\al={J}_\al\circ\dev_* \ (\al=1,2,3)$.
(In fact, noting $H^1(Y;\RR^3)=0$,
there is a map $\lam:Y\ra \RR^3$ for which  
$\displaystyle \dev\bigl((\mbi{t},y)\bigr)=\bigl((\mbi{t}+\lam(y),\Dev(y)\bigr)$
is a gauge transformation from $X=\RR^3\times Y$ onto $\cM$.)
By the uniformization principle, we have the equivariant $qc$-developing pair
$\displaystyle (\rho,\dev):(\Aut_{qc}(X), X)\lra (\Aut_{qc}(\cM),\cM)$
where $\rho$ is the holonomy homomorphism.
Noting $\dev^*\om_0=\eta$, Corollary \ref{abeliabcase}
implies the holonomy homomorphism
reduces to $\displaystyle \rho:\Psh_{qc}(X,\eta)\ra
\Psh_{qc}(\cM,\om_0)\,(=\cM\rtimes (\Sp(n)\cdot \Sp(1)))$. 
As $g_0=\sum_{\al=1}^3 \om_\al\cdot \om_\al+d\om_1\circ J_1$,
$\displaystyle \dev :(X,g)\ra (\cM,g_0)$ turns to a $qc$-local isometry, that is
$\dev^*g_0=g$.
If $M=X/\Gamma$ is compact, then $\dev:X\ra \cM$
is an isometry and hence $M\cong \cM/\rho(\Gamma)$
where $\rho(\Gamma)$ is discrete uniform in $\cM\rtimes (\Sp(n)\cdot \Sp(1))$.
By the generalized Bieberbach theorem, 
$M$ is finitely covered by a quaternionic Heisenberg nilmanifold
$\cM/\Delta$ $(\Delta\leq \cM)$.
\end{proof}

\begin{pro}\label{conthomoqc}
Let $\displaystyle (X,g,\eta,\{J_\al\}_{\al=1}^{3})$ be
a positive definite $4n+3\,(\geq 11)$-dimensional homogeneous $qc$-manifold $G/H$ 
with $qck$-group $\RR^3$.
If $G\leq \Psh_{qc}(X,\eta)$ is a unimodular group,
then $X$ is $qc$-isometric to $\cM$.
\end{pro}

\begin{proof}
As $\RR^3$ is normal in $G$, the quotient manifold
$Y=X/\RR^3$ is homogeneous by the unimodular group $G/\RR^3$.
Then the homogeneous hyperK\"ahler manifold
$\displaystyle (Y,\Omega,\{\hat J_\al\}_{\al=1}^{3})$ 
has a positive definite Hermitian metric $\hat g$
for each $\hat J_\al$. If a direct factor of $\displaystyle (Y,\hat g,\hat J_\al)$
is not flat, then it is a homogeneous K\"ahler manifold of non-compact type
whose Hermitian form has negative Ricci tensor (\cf \cite{HA}, \cite{KN}).
Since $Y$ is homogeneous hyperK\"ahler, applying the proof of Proposition \ref{Dristen},
there is a quaternionic isometry $\Dev: Y\ra\HH^n$
which lifts to a $qc$-isometry $\dev$ of $X$ onto $\cM$.
As $\dev$ is equivariant, there is an isomorphism
 $\displaystyle \rho: \Psh_{qc}(X,\eta)\ra \cM\rtimes(\Sp(n)\cdot \Sp(1))$. 
Since $G$ is transitive on $X$,
$\rho$ induces a $qc$-isometry $\hat\rho: X=G/H\lra \cM$.
\end{proof}

\begin{pro}\label{homoqc1}
Let $M$ be a $4n+3\,(\geq 11)$-dimensional
positive definite locally homogeneous \emph{closed aspherical}  
$qc$-manifold $X/\Gamma$ $(\Gamma\leq \Psh_{qc}(X))$
with a $qck$-distribution $\displaystyle \hat V$.
Then $M$ is $qc$-isometric to 
$\cM/\rho(\Gamma)$ where $\rho(\Gamma)\leq \cM\rtimes \Sp(n)$. Moreover,
$\displaystyle \hat V$ generates $T^3$.
\end{pro}

\begin{proof} 
Since $\Psh_{qc}(M)$ is compact, 
$\hat V$ generates a compact $qck$-subgroup 
$K\leq \Psh_{qc}(M)$.
By the result of
\cite[Theorem 2.4.2, Corollary 3.1.12]{LR},
$K$ is isomorphic to $T^3$. 
Moreover, the orbit map  $\iota(t)=t\cdot x$ $({}^\forall\,t\in T^3)$
at any point $x\in M$ 
induces an injective homomorphism $\iota_*: \ZZ^3=\pi_1(T^3)\ra\,
 C(\Gamma)\leq \Gamma=\pi_1(M)$
where $C(\Gamma)$ is the center of $\Gamma$. 
$X$ inherits a covering group action $\displaystyle \ZZ^3\ra \tilde K\lra K=T^3$
and hence $\tilde K =\RR^3$ which
is a $qck$-subgroup of $\Psh_{qc}(X)$ since so is $K\leq
\Psh_{qc}(M)^0$. 
As $\Gamma\leq \Psh_{qc}(X)$ is cocompact, $\Psh_{qc}(X)$ is unimodular.
Applying Proposition \ref{conthomoqc},
there is a $qc$-isometry $\dev:X\ra \cM$ which induces
a $qc$-isometry of $X/\Gamma$ onto $\cM/\rho(\Gamma)$.
Since $K\leq \Psh_{qc}(M)^0$, 
$\tilde K=\RR^3\leq C_{\Psh_{qc}(X)}(\Gamma)$. Thus
$\rho(\Gamma)$ centralizes $\rho(\RR^3)=C(\cM)=\RR^3$.
As the $\Sp(1)$-action conjugates $\RR^3$,
it follows $\displaystyle \rho(\Gamma)\leq\cM\rtimes \Sp(n)$.
\end{proof}

\begin{pro}\label{compactuni}
Let $\displaystyle \RR\ra\, X/\RR^2\stackrel{q_1}\lra\, Y$
be the pseudo-Hermitian 
bundle as in $(2)$ of Theorem $\ref{hcr}$ where 
$(X/\RR^2,\hat\eta_\al,\hat J'_\al)$
is a $4n+1\,(\geq 9)$-dimensional $CR$-manifold 
with Reeb field $\hat\xi_\al$ 
generating $\RR\leq \RR^3$ $(\al=1,2,3)$.
If $X/\RR^2$ is \emph{spherical}, then $X$ is locally $qc$-isometric to $\cM$.
\end{pro}

\begin{proof}Since $X/\RR^2$ is spherical $CR$,
the pseudo-Hermitian bundle 
shows $(Y,g)$ is a
Bochner flat K\"ahler manifold where $\displaystyle g=\Om_\al\circ \sfJ_\al$,
$d\hat\eta_\al=q_1^*\Om_\al$
(\cf \cite{WE}, \cite{KA1}).
By the result of \cite{TA},
any Bochner flat K\"ahler Einstein manifold 
is the space of constant holomorphic sectional curvature.
As $(Y,g)$ is hyperK\"ahler,
$Y$ is locally holomorphically 
isometric to the flat space $\CC^{2n}$.
Applying the proof of Proposition \ref{Dristen},
$X$ is locally $qc$-isometric to $\cM$.
\end{proof}

\end{document}